\documentclass[12pt]{article}

\usepackage{amsmath}
\usepackage{amssymb}
\usepackage{amsfonts}
\usepackage{a4}
\usepackage{graphicx}

\newtheorem{theorem}{Theorem}[section]
\newtheorem{lemma}[theorem]{Lemma}

\newtheorem{corollary}[theorem]{Corollary}

\newcommand{\Proof}{\par\noindent{\em Proof. }}
\newcommand{\eop}{\nopagebreak\hspace*{\fill}$\Box$\smallskip}

\newcommand{\N}{\Bbb N}
\newcommand{\Z}{\Bbb Z}
\newcommand{\R}{\Bbb R}

\def\Id{\mathbf{Id}}
\def\id{\mathbf{id}}
\def\calL{\mathcal{L}}
\def\eps{\varepsilon}
\def\vv{\mathbf{v}}
\def\trace{\operatorname{trace}}
\def\tW{\tilde{W}}
\def\e{\mathbf{e}}

\def\weakly{\rightharpoonup}

\def\dist{\operatorname{dist}}

\def\XXint#1#2#3{{\setbox0=\hbox{$#1{#2#3}{\int}$}
     \vcenter{\hbox{$#2#3$}}\kern-.5\wd0}}

\usepackage{color}

\begin{document}

\begin{center}
\begin{Large}
{\bf {An atomistic-to-continuum analysis 
of crystal cleavage in a 
two-dimensional model problem}}
\end{Large}
\end{center}

\begin{center}
\begin{large}
Manuel Friedrich\footnote{ Universit{\"a}t Augsburg, Institut f{\"u}r Mathematik, 
Universit{\"a}tsstr.\ 14, 86159 Augsburg, Germany. {\tt manuel.friedrich@math.uni-augsburg.de}}
and Bernd Schmidt\footnote{Universit{\"a}t Augsburg, Institut f{\"u}r Mathematik, 
Universit{\"a}tsstr.\ 14, 86159 Augsburg, Germany. {\tt bernd.schmidt@math.uni-augsburg.de}}\\
\end{large}
\end{center}

\begin{center}
\today
\end{center}
\bigskip

\begin{abstract}
A two-dimensional atomic mass spring system is investigated for critical fracture loads and its crack path geometry. We rigorously prove that, in the discrete-to-continuum limit, the minimal energy of a crystal under uniaxial tension leads to a universal cleavage law and energy minimizers are either homogeneous elastic deformations or configurations that are completely cracked and do not store elastic energy. Beyond critical loading, the specimen generically cleaves along a unique optimal crystallographic hyperplane. For specific symmetric crystal orientations, however, cleavage might fail. In this case a complete characterization of possible limiting crack geometries is obtained. 
\end{abstract}
\bigskip

\begin{small}

\end{small}

\tableofcontents

\section{Introduction}

The behavior of brittle materials is of great interest in applications as well as from a theoretical point of view. Such materials show an elastic response to very small displacements and develop cracks already at moderately large strains. In particular, there is typically no plastic regime in between the restorable elastic deformations and complete failure due to fracture. Major challenges in the experimental sciences and theoretical studies are to identify critical loads at which such a body fails and to determine the geometry of crack paths that occur in the fractured regime. 

In variational fracture mechanics displacements and crack paths are determined from an energy minimization principle. Following the pioneering work of Griffith \cite{Griffith:1921}, Francfort and Marigo \cite{Francfort-Marigo:1998} have introduced an energy functional comprising elastic bulk contributions for the unfractured regions of the body and surface terms that assign energy contributions on the crack paths comparable to the size of the crack of codimension one. Subsequently these models have been investigated and extended in various directions. Among the vast body of literature we only mention the work of Dal Maso and Toader \cite{DalMaso-Toader:2002}; Francfort and Larsen \cite{Francfort-Larsen:2003}; Dal Maso, Francfort and Toader \cite{DalMaso-Francfort-Toader:2005}. Determining energy minimizers of such functionals leads to solving a free discontinuity problem in the language of Ambrosio and De Giorgi \cite{DeGiorgi-Ambrosio:1988} as the crack path, i.e., the set of discontinuity of the diplacement field is not pre-assigned but has to be found as a solution to the variational problem. In particular, these models also lead to efficient numerical approximation schemes, cf., e.g., \cite{Ambrosio-Tortorelli:1992,Bourdin-Francfort-Marigo:2000,Negri:2003,Negri:2005,Schmidt-Fraternali-Ortiz:2009}.

Due to the crystalline structure of matter, under tensile boundary loads fracture typically occurs in the form of cleavage along crystallographic hyperplanes of the atomic lattice. On the continuum side such behavior can be modelled by anisotropic surface terms which are locally minimized for these crack geometries, see e.g.\ \cite{Alicandro-Focardi-Gelli:2000,Buttazzo:1995,Negri:2003}. A discrete model has been investiged by Braides, Lew and Ortiz \cite{Braides-Lew-Ortiz:06}, who assume that fracture can only occur in these directions and then calculate a limiting continuum energy: a cleavage law. This assumption leads to an effective one-dimenional problem which is much easier to analyze. Indeed in the one-dimensional setting, where discrete models describe the behavior of atom chains, a number of results have appeared rather recently on the literature, including \cite{Braides:2000,Braides-Cicalese:2007,Braides-DalMaso-Garroni:1999,Braides-Gelli:2002-1}. While by now for many atomistic models the passage to effective continuum models is well understood in the regime of purely elastic interactions, see \cite{Blanc-LeBris-Lions:2007,Braides-Solci-Vitali:07,Schmidt:2009}, not much is known on discrete-to-continuum limits for models allowing for fracture in more than one dimension. The farthest reaching results in that direction seem to be due to Braides and Gelli \cite{Braides-Gelli:2002-2}, who prove $\Gamma$-convergence results for scalar valued free discontinuity problems. 

However, all these ansatzes fall short of rigorous arguments that indeed in more than one dimension, if fracture occurs at all, then it is energetically favorable to cleave the specimen along particular crystallographic hyperplanes. The main goal of this paper is to provide a rigorous and rather complete study of a two dimensional model problem for validity and failure of crystal cleavage with vector valued deformations. We assume that the body is a rectangular strip subject to uniaxial tensile boundary conditions. Such boundary conditions correspond to one of the basic experiments in determining e.g.\ the Poisson ratio. Our focus on tensile boundary conditions is naturally motivated by our main goal of analyzing cleavage behavior. The atoms in their reference configuration shall be given by the portion of a triangular lattice in that strip that interact via next neighbor Lennard-Jones type potentials. This model seems to be the simplest model problem which (1) is frame indifferent in its vector-valued arguments in more than one dimension, (2) gives rise to non-degenerate elastic bulk terms and (3) leads to surface contributions sensitive to the crack geometry with competing crystallographic hyperplanes. Moreover, two-dimensional lattice surfaces naturally appear in the analysis of thin structures. In particular we will also discuss consequences of our analysis on the stability of brittle nanotubes under interior expansive pressure.

Indeed we will prove that under uniaxial tension in the continuum limit the energy satisfies a particular cleavage law with quadratic response to small boundary displacements followed by a sharp constant cut-off beyond some critical value. Moreover, we will see that any sequence of minimizers converges (up to subsequences) to a homogeneous continuum deformation for subcritical boundary values, while it converges to a continuum deformation which is cracked completely and does not store elastic energy in the supercritical case. We show that in the generic case cleavage occurs along a unique crystallygraphic line, whereas for specific symmetric orientations of the crystal cleavage might fail. Nevertheless, also in these special cases we obtain a complete characterization of all possible limiting crack geometries. The model under investigation leads, in particular, to configurations respecting the Poisson effect, which would not be possible in scalar models. These results justify rigorously the aforementioned assumptions in the derivation of cleavage laws as, e.g., in \cite{Braides-Lew-Ortiz:06}. 

Of course the derivation of continuum theories for brittle materials from atomistic interactions remains a challenging open problem in its full generality. To the best of our knowledge, however, we believe that our analysis provides a first vectorial result in more than one dimension. For our model specimen we can completely classify its response to lateral stretching according the strength of stretching and all possible lattice orientations and also prove strong discrete-to-continuum convergence results in the various regimes. Even though the uniaxial tension test is a natural set-up for investigating cleavage phenomena, it would be desirable to also incorporate more genaral boundary conditions. This appears to be a rather difficult task. Even in the case of uniaxial compressions our model would predict unphysical behavior and one would have to make further modelling assumptions corresponding to an assumption in a one-dimensional situation where atoms are not allowed to pass each other. An even more challenging open problem is to identify limiting continuum configurations and energies in the same energy regime which are not necessarily asymptotically energy minimizing. We only have partial results in this direction so far, cf.\ \cite{FriedrichSchmidt-Gamma}. Here a limiting variational principle is derived under additional suitable smallness assumptions on the discrete displacement, whose minimizers are exactly the minimizers found in the present analysis. On the other hand, we believe that our techniques do allow for generalizations in several other directions, in particular to more general Bravais lattices, more than two dimensions and longer range interaction potentials. We defer such investigations to subsequent work.  

The paper is organized as follows. We first introduce our discrete model and state our main results in Section \ref{sec:model-and-main-results}. Here we already discuss different scalings of the boundary data and find the interesting regime where both energy contributions, the elastic and the crack energy, are of the same order. 

In Section \ref{sec:cell-energy} we collect some elementary properties of the cell energy. In particular, we introduce a lower-bound comparison energy, called \textit{reduced energy}, providing the optimal cell energy in dependence of the cell expansion in the space direction where tensile boundary conditions are imposed. 

Section \ref{sec:limiting-energy} is devoted to the derivation of cleavage laws for the limiting minimal energy. Using an elementary slicing argument we reduce the problem to one-dimensional segments and show that the limiting energy has a universal form independent of the interatomic potential. Our result is similar to the effective one-dimesional law discussed in \cite{Braides-Lew-Ortiz:06}. We obtain that the crack energy is anisotropic and depends explicitly on the lattice orientation. While such a law can be proved with rather elementary methods, more subtle arguments have to be used in order to obtain finer estimates on the limiting energy. We will derive higher order terms for the discrete minimal energies and see that in the limiting behavior of these terms also anisotropic contributions occur in the elastic regime. Moreover, our proof illustrates the typical behavior of brittle materials already seen in the continuum cleavage law also in a discrete framework: There is essentially no plastic regime besides the elastic and the crack regime. More precisely, we see that for almost minimizers the deformation is either $\sqrt{\eps}$-close to the identity mapping (representing elastic response) or springs between adjacent atoms are elongated by a factor scaling like $\frac{1}{\sqrt{\eps}}$ (leading to fracture in the limit description), where $\eps$ denotes the typical interatomic distance. In particular, here we can already see that homogeneous deformations or cleavage along specific lines are asymptotically optimal and that for specific symmetric crystal orientation the crack geometry may become much more complicated. 

In Section \ref{sec:limiting-configurations} we proceed to show that, under appropriate assumptions, in terms of suitably rescaled displacement fields indeed all discrete energy minimizers converge strongly to such continuum deformations. The main difficulty to establish a compactness result in the vector valued case arises from the geometric nonlinearity induced by the frame indifference of the model. In a purely elastic regime such a compactness for configurations with deformation gradient near the orientation preserving rigid motions is due to DalMaso, Negri and Percivale \cite{DalMasoNegriPercivale:02} for continuum deformations and the second author \cite{Schmidt:2009} for discrete interactions by applying the geometric rigidity theorem of Friesecke, James and M{\"u}ller \cite{FrieseckeJamesMueller:02}. Such a rigidity estimate, however, is not adequate for the framework of brittle materials as we do not have $p$-growth conditions for the atomic interactions. To derive a rigidity result for our problem we overcome this difficulty by providing a fine characterization of the crack, i.e. of the number and position of largely elongated springs. In the subcritical case the contribution of such springs is abitrarily small such that the purely elastic theory applies. In the generic case, for supercritical boundary values largely deformed springs lie in a small stripe in direction of the optimal cristallographic line in such a way that the two components on the right and on the left of the stripe essentially behave elastically. We then also establish a strong convergence result.

\section{The model and main results}\label{sec:model-and-main-results}

\subsection*{The discrete model}

Let ${\cal L}$ denote the rotated triangular lattice 
$$ {\cal L} 
   = R_{\cal L} \begin{pmatrix} 1 & \frac{1}{2} \\ 0 & \frac{\sqrt{3}}{2} \end{pmatrix} \Z^2 
   = \{ \lambda_1 \vv_1 + \lambda_2 \vv_2 : \lambda_1, \lambda_2 \in \Z \}, $$ 
where $R_{\cal L} \in SO(2)$ is some rotation and $\vv_1$, $\vv_2$ are the lattice vectors $\vv_1 = R_{\cal L} \mathbf{e}_1$ and $\vv_2 = R_{\cal L}(\frac{1}{2} \mathbf{e} + \frac{\sqrt{3}}{2} \mathbf{e}_2)$, respectively. We collect the basic lattice vectors in the set ${\cal V} = \left\{\vv_1,\vv_2,\vv_2 - \vv_1\right\}$. The region $\Omega = (0, l) \times (0, 1) \subset \R^2$, $l > 0$, is considered the macroscopic region occupied by the body under investigation. In the reference configuration the positions of the specimen's atoms are given by the points of the scaled lattice $\eps\calL$ that lie within $\Omega$. Here $\eps$ is a small parameter defining the length scale of the typical interatomic distances. 

The deformations of our system are mappings $y : \eps \calL \cap \Omega \to \R^2$. The energy associated to such a deformation $y$ is assumed to be given by nearest neighbor interactions as 
\begin{align}\label{eq:Energy}
  E_{\eps}(y) 
  = \frac{1}{2} \sum_{x,x' \in \eps {\cal L} \cap \Omega \atop |x-x'| = \eps} W \left( \frac{|y(x) - y(x')|}{\eps} \right).  
\end{align} 
Note that the scaling factor $\frac{1}{\eps}$ in the argument of $W$ takes account of the scaling of the interatomic distances with $\eps$. The pair interaction potential $W:[0,\infty) \to [0, \infty]$ is supposed to be of `Lennard-Jones-type': 
\begin{itemize}\label{W-assumptions}
\item[(i)] $W \ge 0$ and $W(r) = 0$ if and only if $r = 1$. 
\item[(ii)] $W$ is continuous on $[0, \infty)$ and $C^2$ in a neighborhood of $1$ with $\alpha := W''(1) > 0$. 
\item[(iii)] $\lim_{r \to \infty} W(r) = \beta$. 
\end{itemize} 
In order to obtain fine estimates on limiting energies and configurations we will also consider the following stronger versions of hypotheses (ii) and (iii): 
\begin{itemize}
\item[(ii')] $W$ is continuous on $[0, \infty)$ and $C^4$ in a neighborhood of $1$ with $\alpha := W''(1) > 0$ and arbitrary $\alpha' := W'''(1)$.
\item[(iii')] $W(r) = \beta + O(r^{-2})$ as $r \rightarrow \infty$,  
\end{itemize} 
which is still satisfied, e.g., by the classical Lennard-Jones potential. 
 
In order to analyze the passage to the limit as $\eps \to 0$ it will be useful to interpolate and rewrite the energy as an integral functional. Let ${\cal C}_{\eps}$ be the set of equilateral triangles $\triangle \subset \Omega$ of sidelength $\eps$ with vertices in $\eps {\cal L}$ and define $\Omega_{\eps} = \bigcup_{\triangle \in {\cal C}_{\eps}} \triangle$. By $\tilde{y} : \Omega_{\eps} \to \R^2$ we denote the interpolation of $y$, which is affine on each $\triangle \in {\cal C}$. The derivative of $\tilde{y}$ is denoted by $\nabla \tilde{y}$, whereas we write $(y)_{\triangle}$ for the (constant) value of the derivative on a triangle $\triangle \in {\cal C}_{\eps}$. Then \eqref{eq:Energy} can be rewritten as 
\begin{align}\label{eq:E-integral}
\begin{split}
  E_{\eps}(y) 
  &= \sum_{\triangle \in {\cal C}_{\eps}} W_{\triangle} ((\tilde{y})_{\triangle}) 
     + E_{\eps}^{\rm boundary}(y) \\ 
  &= \frac{4}{\sqrt{3}\eps^2} \int_{\Omega_{\eps}} W_{\triangle} (\nabla \tilde{y}) \, dx 
     + E_{\eps}^{\rm boundary}(y), 
\end{split} \end{align}
where 
\begin{align}\label{eq:W-triangle}
  W_{\triangle}(F) 
  &= \frac{1}{2} \Big( W(|F \vv_1|) + W(|F \vv_2|) + W(|F (\vv_2 - \vv_1)|) \Big).     
\end{align} 
(Note that $|\triangle| = \sqrt{3}\eps^2/4$.) Here the boundary term is the sum of pair interaction energies $\frac{1}{4} W(\frac{|y(x) - y(x')|}{\eps})$ over nearest neighbor pairs which form the side of none or only one triangle in ${\cal C}_{\eps}$.

\subsection*{Boundary values and scaling}

We are interested in the behavior of the system when applying tensile boundary conditions, say in $\mathbf{e}_1$-direction. In particular, we would like to investigate when and how the body breaks, i.e., 
\begin{itemize}
\item[(1)] at which value of the boundary displacement energetic minimizers are no longer elastic deformations but exhibit cracks and 
\item[(2)] if indeed it is most favorable for the cracks to separate the body along crystallographic lines. 
\end{itemize} 
In order to avoid geometric artefacts, we will therefore assume that $l > \frac{1}{\sqrt{3}}$, so that it is possible for the body to completely break apart along lines parallel to $\R\vv_1$, $\R\vv_2$ or $\R(\vv_2 - \vv_1)$ not passing through the left or right boundaries. 

Due to the discreteness of the underlying atomic lattice we have to impose the boundary conditions of uniaxial extension in small neighborhoods of $\{0\} \times (0,1)$ and $\{l\} \times (0, 1)$, respectively, as otherwise cracks near the boundary may become energetically more favorable. For $a_{\eps} > 0$ we set 
\begin{align*} 
\begin{split}
  {\cal A}(a_{\eps}) 
  &= \big\{ y = (y_1, y_2) : \eps {\cal L} \cap \Omega \to \R^2 : \\ 
  &\qquad \qquad \qquad \qquad  
     y_1(x) = (1 + a_{\eps}) x_1 \text{ for } x_1 \le \eps \text{ and } x_1 \ge l - \eps \big\}. 
     \end{split}
\end{align*}
In the special case $\phi = 0$ we will in addition assume that there is an upper bound $R_0$ on the elongation of every atomic bond in a small $\psi(\eps)$-neighborhood of the lateral boundaries of width $\psi(\eps) > 0$ with $\eps \ll \psi(\eps) \ll 1$: 
\begin{align} \label{eq: bv-phi-0} 
  |y(x) - y(x')| \le R_0 \eps 
  \text{ if } x_1, x_1' \le \psi(\eps) \text{ or } x_1, x_1' \ge l - \psi(\eps). 
\end{align}
Without such an assumption, in the general low energy regime to be considered later, for $\phi = 0$ the boundary values are not strong enough to prevent the specimen from breaking on the boundary into a large amount of completely separated components, rendering the system too sensitive to unphysical boundary effects. 

Note that, apart from choosing $\psi$, there is some arbitrariness in this implementation of boundary values as one might, e.g., equally well ask that 
\begin{align*}
  y_1(x) = x_1 \text{ for } x_1 \le \eps ~\text{ and }~ 
  y_1(x) = x_1 + a_{\eps} l \text{ for } x_1 \ge l - \eps. 
\end{align*}
Such different choices will, however, not change the results of the analysis. 

Also note that there is no assumption on the second component of the boundary displacement, i.e., the atoms may `slide along the boundary lines'. Besides describing a basic experiment on elastic bodies, this assumption allows for a direct application of our results to the stability analysis of nanotubes:

If the rotation $R_{\cal L}$ and the length $l$ are such that for a sequence $\eps_k \to 0$ the translated lattice $\eps_k {\cal L} + (l, 0)$ concides with the original lattice $\eps_k{\cal L}$, we may view the system as an atomistic nanotube with macroscopic region $\frac{l}{2\pi}S^1 \times (0, 1)$. (Note that for small $\eps_k$ the bending energy contributions when rolling up $(0, l) \times (0, 1)$ into a cylinder are negligible as this mapping is an isometric immersion and thus infinitesimally rigid.) Imposing periodic boundary conditions, for arbitrary $l > 0$ our system then models deformations of a nanotube subject to expansion of the diameter.  

There are two obvious choices for deformations satisfying the boundary conditions: The homogeneous elastic deformation $y^{\rm el}(x) = (1 + a_{\eps}) x$ and a cracked body deformation $y^{\rm cr}$, which, up to a boundary layer of negligible energy, is the identity to the left and a translation by $a_{\eps} l \e_1$ to the right of some segment (or curve) passing through $\Omega$ that connects a point on the lower boundary $(\eps,l- \eps) \times \{0\}$ and a point on the upper boundary $(\eps,l- \eps) \times \{1\}$. It is not hard to see that 
$$ E_{\eps}(y^{\rm el}) \sim \eps^{-2} W(1 + a_{\eps}), \quad \quad 
   E_{\eps}(y^{\rm cr}) \sim \eps^{-1}. $$ 
In particular, we are interested in the most interesting regime where both of these energy values are of the same order, i.e., $a_{\eps}$ is small and 
$$ \eps^{-2} a_{\eps}^2 \sim \eps^{-2} W(1 + a_{\eps}) \sim \eps^{-1} 
   \quad \implies \quad 
   a_{\eps} \sim \sqrt{\eps}. $$
In order to obtain finite and nontrivial energies in the limit $\eps \to 0$, we accordingly rescale $E_{\eps}$ to ${\cal E}_{\eps} := \eps E_{\eps}$. 

Conceivable alternative implementations of the boundary conditions as alluded to above will then result in energy changes of order $O(\eps)$. We will account for all such possibilities by characterizing not only energy minimizing configurations, but more generally all configurations which are energy minimizing up to an error term of order $O(\eps)$.

\subsection*{Limiting minimal energy and cleavage laws}

We begin our analysis with an elementary argument which yields the limiting minmal energy as $\eps \to 0$ when $a_{\eps}/\sqrt{\eps} \to a \in [0, \infty]$. We first establish a lower bound for this energy by considering slices of the form $(0,l) \times \{x_2\}$ for $x_2 \in (0,1)$ and using the {\em reduced energy} $\tW$ defined by  
\begin{equation}\label{eq: reduced energy} 
\tW(r) = \inf\{ W_{\triangle}(F) : \e_1^T F \e_1 = r\}. 
\end{equation}
In a second step we show that this bound is attained. In particular, it turns out that the limiting minimal energy is given by elastic deformations up to some critical value $a_{\rm crit}$ of the boundary displacements and by cleavage along a specific crystallographic line beyond $a_{\rm crit}$. 

Let $\gamma = \max\{|\vv_1 \cdot \e_2|, |\vv_2 \cdot \e_2|, |(\vv_2 - \vv_1) \cdot \e_2|\}$ and $\vv_{\gamma} \in {\cal V}$ such that $\gamma = |\vv_{\gamma} \cdot \e_2|$. We note that $\gamma$ takes values in $[\sqrt{3}/2,1]$ and that $\vv_{\gamma}$ is unique if $\gamma > \sqrt{3}/2$. 
\begin{theorem}\label{theo:limiting-energy} 
Suppose $a_{\eps}/\sqrt{\eps} \to a \in [0, \infty]$. The limiting minimal energy is given by 
\begin{equation}\label{eq:limiting-energy} 
  \lim_{\eps \to 0} \inf \left\{ {\cal E}_{\eps}(y) : y \in {\cal A}(a_{\eps}) \right\} 
  = \min \left\{\frac{ \alpha l}{\sqrt{3}} a^2, \frac{2 \beta}{\gamma} \right\}. 
\end{equation}
\end{theorem}

As already motivated above, only one of the regimes is energetically favorable if $a \in \{0, \infty\}$. In the interesting case $a \in (0, \infty)$ we indeed will see that in terms of the critical boundary displacement 
$$ a_{\rm crit} = \sqrt{\frac{2 \sqrt{3} \beta}{ \alpha \gamma l}} $$ 
the limit is attained for homogeneously deformed configurations if $a \le a_{\rm crit}$ and for configurations cracked along lines parallel to $\R \vv_{\gamma}$, if $a \ge a_{\rm crit}$. In the special case that $\vv_{\gamma}$ is not unique the limit is also attained if the crack takes a serrated course parallel to $\R (\frac{1}{2},\frac{\sqrt{3}}{2})^{T}$ or  $\R(-\frac{1}{2},\frac{\sqrt{3}}{2})^{T}$. 

For the sake of simplicity we specialize to sequences $a_{\eps} = \sqrt{\eps} a$. Without loss of generality we assume that $R_{\cal L} = \begin{footnotesize} \begin{pmatrix} \cos \phi & - \sin \phi \\ \sin \phi & \cos \phi \end{pmatrix} \end{footnotesize}$ for $\phi \in [0, \frac{\pi}{3})$, so that $\gamma = \sin( \phi + \frac{\pi}{3}) = \vv_{\gamma} \cdot \e_2$ and $\vv_{\gamma}$ is unique iff $\phi \ne 0$. If the assumptions (ii') and (iii') on $W$ hold, we have the following sharp estimate on the discrete minimal energies up to error terms of the order of surface contributions. 
\begin{theorem}\label{theo: discrete energy}
For $\eps$ small the discrete minimal energy is given by 
$$ \inf{\cal E}_{\eps} = \min \left\{ \frac{ \alpha l}{\sqrt{3}} a^2 + \frac{[6 \alpha + 7 \alpha' - 2(3 \alpha - \alpha') \cos(6 \phi)] l}{27 \sqrt{3}} \sqrt{\eps} a^3 , \frac{2 \beta}{\gamma} \right\} + O(\eps). $$
\end{theorem}

Thus, while the zeroth order contributions in the elastic regime are isotropic, the higher order contributions as well as the fracture energy  explicitly depend on the lattice orientation angle $\phi$. 

Detailed proofs of these results will be given in Section \ref{sec:limiting-energy}.

\subsection*{Limiting minimal configurations}

Our analysis of the limiting minimal energy so far showed that, depending on the boundary data, homogeneous deformations or completely cracked configurations are energy minimizing in the limit $\eps \to 0$. However, it falls short of showing that in fact these configurations are the only possibilities to obtain asymptotically optimal energies. Indeed, if $\vv_{\gamma}$ is not unique, then we have already seen that the crack path can become geometrically much more complicated. Our next result shows that energy minimizing configurations converge to a homogeneous continuum deformation for subcritical boundary values, while  in the supercritical case they converge to a continuum deformation which is completely cracked and does not store elastic energy. If $\phi \ne 0$ and hence $\vv_{\gamma}$ is unique, the crack path follows the optimal crystallographic line. For $\phi = 0$ such a cleavage behavior fails in general. Nevertheless, we obtain an explicit characterization of all possible limiting crack shapes in this case as well in terms of Lipschitz curves whose tangent vector lies a.e.\ in the cone generated by $(-\frac{1}{2},\frac{\sqrt{3}}{2})$ and $(\frac{1}{2},\frac{\sqrt{3}}{2})$.

The basic idea behind our reasoning will be to `count' the number of `broken' springs, i.e. the springs intersected transversally by the crack path. We see that the springs broken by a crack line $(p,0) + \R\vv_{\gamma}$ do not overlap in the projection onto the $x_2$-axis and the length of the projection of two adjacent broken springs equals $\eps \gamma$. This leads to a fracture energy of approximately $\frac{2 \beta}{\gamma}$. If in the generic case $\phi \neq 0$ we assume that the cleavage is not parallel to $\R\vv_{\gamma}$ we conclude that some springs in $\vv_{\gamma}$ direction must be broken, too. If we consider the adjacent triangles of such a spring and their neighbors we find that the projection onto the $x_2$-axis of broken springs overlap. A careful analysis of this phenomenon then shows that every broken spring in $\vv_{\gamma}$ direction `costs' an additional energy of $\approx 2 \eps \beta \frac{P(\gamma)}{\gamma}$, where $P(\gamma)$ is the geometrical factor 
\begin{align}\label{eq:P-gamma}
  P(\gamma) 
  = \frac{1}{2} \left(1- \sqrt{3} \frac{\sqrt{1 - \gamma^2}}{\gamma}\right). 
\end{align}
(Note that $P(\gamma) = 0 \Leftrightarrow \gamma = \frac{\sqrt{3}}{2} \Leftrightarrow \phi = 0$ in accordance to the above considerations.) For the special case $\phi = 0$ we provide a similar counting argument. 

In order to give a precise meaning to the convergence of discrete to continuum deformations, to each discrete deformation $y : \eps \calL \to \R^2$ we assign -- as mentioned above -- the affine interpolation $\tilde{y}$ on each triangle $\triangle \in {\cal C}_{\eps}$. Accordingly, to the rescaled discrete displacements $u: \eps \calL \to \R^2$ with $y = \id + \sqrt{\eps} u$ ($\id$ denoting the identity mapping $\id(x) = x$) we define $\tilde{u}$ to be its affine interpolation on each triangle $\triangle \in {\cal C}_{\eps}$. 

In the cracked regime we may of course only hope for a unique limiting deformation up to translation of the crack path. However, without an additional mild extra assumption on the admissible discrete configurations or their energy even this cannot hold true, as apart from the crack, parts of the specimen could flip their orientation and fold onto other parts on the body at zero energy. In order to avoid such unphysical behavior we add a frame indifferent penalty term $\chi \ge 0$ to $W_{\triangle}$ with $\chi \geq c_{\chi} > 0$ in a neighborhood of $O(2) \setminus SO(2)$ and $\chi \equiv 0$ in a neighborhood of $SO(2)$ and $\infty$, which in particular does not change the energy response in the linear elastic and in the fracture regime:
\begin{align}\label{eq:W-trinagle-chi}
  W_{\triangle, \chi}(F) = W_{\triangle}(F) + \chi(F). 
\end{align}
For instance, an admissible choice for $\chi$ is the local orientation preserving condition in the elastic regime 
$$\chi(F) = \begin{cases} 0,  & \text{if } \det(F) > 0 \text{ or } |F| > R, \\
                          \infty, & \text{if } \det(F) \leq 0 \text{ and } |F| \leq R, \end{cases}$$
for some threshold $R \gg 1$. (Also $1 \le R = R(\eps) \ll \frac{1}{\sqrt{\eps}}$ would be admissible.) We remark that such an infinitely strong penalization of deformation gradients with non-positive determinant is widely used in the elastic models. Allowing for fracture, however, a penalization of orientation reversion between different cracked parts of the body is no longer physically justifiable, whence we set $\chi = 0$ for very large deformation gradients.
We set 
\begin{align*}
  {\cal E}^{\chi}_{\eps} (y) 
  = \frac{4}{\sqrt{3}\eps} \int_{\Omega_{\eps}} W_{\triangle, \chi} (\nabla \tilde{y}) \, dx + \eps E^{\rm boundary}_{\eps}(y), 
\end{align*}
for $u \in {\cal A}(a_{\eps})$. More generally than a sequence of minimizers we will consider sequences $(y_{\eps})$ of almost minimizers that satisfy 
\begin{align}\label{eq:almost-min}
  {\cal E}^{\chi}_{\eps} (y_{\eps}) 
  = \inf\{ {\cal E}^{\chi}_{\eps}(y) : y \in {\cal A}(a_{\eps}) \} + O(\eps). 
\end{align}
For those deformations we will show in Section \ref{sec:limiting-configurations}:
\begin{theorem}\label{theo: conv-of-minimizers} 
Assume that $W$ satisfies (i), (ii') and (iii'). Let $a_{\eps} = \sqrt{\eps} a$, $a \ne a_{\rm crit}$ and suppose $(y_{\eps})$ satisfies \eqref{eq:almost-min}. Let $u_{\eps}$ such that $y_{\eps} = \id + \sqrt{\eps} u_{\eps}$. Then there exist $\bar{u}_{\eps} : \Omega \to \R^2$ with $|\{x \in \Omega_{\eps}: \bar{u}_{\eps}(x) \ne \tilde{u}_{\eps}(x)\}| = O(\eps)$ such that: 
\begin{itemize} 
\item[(i)] If $a < a_{\rm crit}$, then there is a sequence $s_{\eps} \in \R$ such that 
$$ \| \bar{u}_{\eps} - (0,s_{\eps}) - F^a \cdot \|_{H^1(\Omega)} \to 0, $$ 
where $F^a = \begin{pmatrix} a & 0  \\ 0 & -\frac{a}{3} \end{pmatrix}$.  
\item[(ii)] If $a > a_{\rm crit}$ and $\phi \ne 0$, then there exist sequences $p_{\eps} \in (0, l)$, $s_{\eps}, t_{\eps} \in \R$ such that $(p_{\eps},0) + \R\vv_{\gamma}$ intersects both the segments $(0, l) \times \{0\}$ and $(0, l) \times \{1\}$ and, for the parts  to the left and right of $(p_{\eps},0) + \R\vv_{\gamma}$
\begin{align*}
  \Omega^{(1)} 
  &:= \left\{x \in \Omega : 0 < x_1 < p_{\eps} + (\vv_{\gamma} \cdot e_1) x_2  \right\} \text{ and } \\ 
  \Omega^{(2)} 
  &:= \left\{x \in \Omega : p_{\eps} + (\vv_{\gamma} \cdot e_1) x_2 < x_1 < l  \right\}, 
\end{align*}
respectively, we have 
$$ \| \bar{u}_{\eps} - (0, s_{\eps}) \|_{H^1(\Omega^{(1)})} 
   + \| \bar{u}_{\eps} - (al, t_{\eps}) \|_{H^1(\Omega^{(2)})} 
   \to 0. $$ 
   \item[(iii)] If $a > a_{\rm crit}$ and $\phi = 0$, then there exist sequences of Lipschitz functions $g_{\eps}: (0,1) \to (0, l)$ satisfying $g'_\eps=\pm\frac{1}{\sqrt{3}}$ a.e.\ such that for the parts to the left and right of $\text{graph}(g_\eps)$ 
\begin{align*}
  \Omega^{(1)}[g_\eps] 
  &:= \left\{x \in \Omega : 0 < x_1 < g_{\eps}(x_2) \right\} \text{ and } \\ 
  \Omega^{(2)}[g_\eps] 
  &:= \left\{x \in \Omega : g_{\eps}(x_2) < x_1 < l  \right\}, 
\end{align*}
respectively, we have 
$$ \| \bar{u}_{\eps} - (0,  s_{\eps}) \|_{H^1(\Omega^{(1)}[g_\eps])} 
   + \| \bar{u}_{\eps} - (al, t_{\eps}) \|_{H^1(\Omega^{(2)}[g_\eps])} 
   \to 0, $$ 
for suitable sequences $s_{\eps}, t_{\eps} \in \R$. 
\end{itemize}
\end{theorem} 

As a consequence, we obtain a complete characterization of limiting continuum deformations, when no mass leaks to infinity.

\begin{corollary}\label{cor: limiting-deformations} 
Under the assumptions and with the notation of Theorem \ref{theo: conv-of-minimizers}, if $\sup_\eps \left\|y_\eps\right\|_\infty < \infty$, up to passing to subsequences, $\tilde{u}_{\eps} \to u$ in measure where 
\begin{itemize} 
\item[(i)] if $a < a_{\rm crit}$, $u(x) = F^a x + (0,s)$ for some constant $s \in \R$, 

\item[(ii)] if $a > a_{\rm crit}$ and $\phi \ne 0$, 
$u(x) 
= \begin{cases} (0, s), & \mbox{for } x \mbox{ to the left of } (p,0) + \R\vv_{\gamma}, \\ 
                (al, t), & \mbox{for } x \mbox{ to the right of } (s,0) + \R\vv_{\gamma}, 
\end{cases}$ 
for constants $s, t \in \R$ and $p \in (0,l)$ such that $(p,0) + \R\vv_{\gamma}$ intersects both the segments $(0,l) \times \left\{0\right\}$ and $(0,l) \times \left\{1\right\}$, 

\item[(iii)] if $a > a_{\rm crit}$ and $\phi = 0$, 
$u(x) 
= \begin{cases} (0, s), & \mbox{if } 0 < x_1 < g(x_2), \\ 
                (al, t), & \mbox{if } g(x_2) < x_1 < l, 
\end{cases}$ 

for a Lipschitz function $g : (0,1) \to [0, l]$ with $|g'| \le \frac{1}{\sqrt{3}}$ a.e.\ and constants $s, t \in \R$. 
\end{itemize}
Conversely, for every $u$ as given in the cases (i)-(iii) there is a minimizing sequence $(y_\eps)$ satisfying \eqref{eq:almost-min} and $\tilde{u}_{\eps} \to u$ in measure.
\end{corollary}

We close this introductory chapter emphasizing that all the optimal configurations found in Theorem \ref{theo: conv-of-minimizers} and Corollary \ref{cor: limiting-deformations} by minimizing the energy without a priori assumptions show purely elastic behavior in the subcritical case and complete fracture in the supercritical regime. In particular, the elastic minimizer in (i) shows elongation $a$ in $\e_1$-direction and compression $-\frac{a}{3}$ in the perpendicular $\e_2$-direction, a manifestation of the Poisson effect (with Poisson ratio $\frac{1}{3}$), which cannot be derived in scalar valued models. On the other hand, the crack minimizer in (ii) for $\phi \ne 0$ is broken parallel to $\R \vv_{\gamma}$ which proves that cleavage occurs along crystallographic lines, while we see that cleavage in the symmetric case $\phi = 0$ in general fails.

\section{Elementary properties of the cell energy}\label{sec:cell-energy}

We collect some properties of the cell energy $W_{\triangle}$ and the reduced energy defined in (\ref{eq: reduced energy}) for $W$ satisfying the assumptions (i), (ii) and (iii).

\begin{lemma}\label{lemma:W-triangle-properties}  
$W_{\triangle}$ is 
\begin{itemize}
\item[(i)] frame indifferent: $W_{\triangle}(QF) = W_{\triangle}(F)$ for all $F \in \R^{2 \times 2}$, $Q \in SO(2)$, 
\item[(ii)] non-negative and satisfies $W_{\triangle}(F) = 0$ if and only if $F \in O(2)$ and 
\item[(iii)] $\liminf_{|F| \to \infty} W_{\triangle}(F) = \liminf_{|F| \to \infty} W_{\triangle,\chi}(F) = \beta$. 
\end{itemize}
\end{lemma}

\Proof (i) is clear. For (ii) it suffices to note that $v^T F^T F v = 1$ for three vectors $v$, no two of which are collinear, implies that $F^T F = \Id$. As $\chi$ vanishes near $\infty$, (iii) can be seen by noting that if $|F| \to \infty$, then for at least two vectors $\vv \in  {\cal V}$ one has $|F\vv| \to \infty$. \eop 

We compute the linearization about the identity matrix $\Id$:

\begin{lemma}\label{lemma:W-linearized-properties}
Let $F = \Id + G$ for $G \in \R^{2 \times 2}$. Then for $|G|$ small 
$$ W_{\triangle}(F) = \frac{1}{2} Q(G) + o(|G|^2),$$
where $Q(G) = \frac{3 \alpha}{16} \left( 3 g_{11}^2 + 3 g_{22}^2 + 2 g_{11} g_{22} + 4 \left(\frac{g_{12}+g_{21}}{2} \right)^2 \right)$. 

In particular, $Q(G)$ only depends on the symmetric part $\left(G^{T} + G\right)/2$ of $G$. $Q$ is positive semidefinite and thus convex on $\R^{2 \times 2}$ and positive definite and strictly convex on the subspace $\R^{2 \times 2}_{\rm sym}$ of symmetric matrices. 
\end{lemma}

\Proof Let $\vv \in {\cal V}$ and $G \in \R^{2 \times 2}$ small. We Taylor expand the contributions $W(|F \vv|)$ to the energy $W_{\triangle}$: 
\begin{align*}
  W(| (\Id +  G) \vv|) 
  &= W \left( \sqrt{\langle \vv, (\Id +  G^T)(\Id +  G) \vv \rangle} \right) \\ 
  &= \frac{W''(1)}{2} \left\langle \vv, \frac{G^T + G}{2} \vv  \right\rangle^2 + o(|G|^2). 
\end{align*} 
Now using the elementary identity 
\begin{align}\label{eq:lat-sum}
\begin{split}
  &\langle \vv_1, H \vv_1 \rangle^2 
  + \langle \vv_2, H \vv_2 \rangle^2 
  + \langle (\vv_2 - \vv_1), H (\vv_2 - \vv_1) \rangle^2 \\
  &= \frac{3}{8} \left( 2\trace(H^2) + (\trace H)^2 \right) 
\end{split}
\end{align}
for any symmetric matrix $H \in \R^{2 \times 2}$, we obtain by summing over $\vv \in {\cal V}$ 
\begin{align*}
  W_{\triangle} (F)
  &= \frac{1}{2} \cdot \frac{\alpha}{2} \cdot \frac{3}{8} \cdot 
     \left( 2 \trace\left(\left(\frac{G^T + G}{2}\right)^2\right) + \left(\trace \frac{G^T + G}{2}\right)^2 \right) 
     + o(|G|^2)\\ 
  &= \frac{1}{2}Q(G) + o(|G|^2).
\end{align*} 
As $Q(G) \ge \frac{3 \alpha}{16} (2 g_{11}^2 + 2 g_{22}^2 + (g_{12} + g_{21})^2)$, $Q$ is positive semidefinite on $\R^{2 \times 2}$ and positive definite on $\R^{2 \times 2}_{\rm sym}$. \eop

As a consequence, we have the following properties of the reduced energy $\tW$. 

\begin{lemma}\label{lemma:W-reduced-properties} 
The reduced energy satisfies  
\begin{itemize}
\item[(i)] $\tW(r) = 0 \iff |r| \le 1$. 
\item[(ii)] For $r \ge 1$ one has 
$$ \tW(r) 
   = W_{\triangle} \left( \begin{pmatrix} r & 0 \\ 0 & \frac{4 - r}{3} \end{pmatrix} \right) + o((r - 1)^2) 
   = \frac{ \alpha}{4} (r - 1)^2 + o((r - 1)^2). $$ 
\item[(iii)] $\lim_{|r| \to \infty} \tW(r) = \beta$. 
\end{itemize}
\end{lemma}

\Proof 
(i) If $r \le 1$, then one can choose $Q \in SO(2)$ with $\e_1^T Q \e_1 = r$ and so $0 \le \tW(r) \le W_{\triangle}(Q) = 0$. If $|r| > 1$, then $\tW(r) > 0$ for otherwise there would be a sequence $F_k \in \R^{2 \times 2}$ with $\e_1^T F_k \e_1 = r$ and $W_{\triangle}(F_k) \to 0$. But then ${\rm dist}(F_k, O(2)) \to 0$ by (ii) and (iii) of Lemma \ref{lemma:W-triangle-properties} and thus, up to subsequences, $F_k \to F \in O(2)$ with $\e_1^T F \e_1 = r$, which is impossible. 

(ii) This discussion shows that in fact for any $\delta > 0$ there exists $\eta > 0$ such that $W_{\triangle}(F) > \delta$ whenever ${\rm dist}(F, O(2)) \ge \eta$. Now since $\tW(r) \to 0$ as $r \searrow 1$, we obtain that, for sufficiently small $r > 1$ and $\delta > 0$, any $F$ with $W_{\triangle}(F) < \tW(r) + \delta$ is contained in a small neighborhood of $O(2)$. If in addition $\e_1^T F \e_1 = r$ holds, then in fact, $F$ must be close to $\Id$ or to $P = \begin{footnotesize}\begin{pmatrix} 1 & 0 \\ 0 & -1 \end{pmatrix}\end{footnotesize}$. In particular, by continuity of $W$, the infimum on the right hand side in the definition of $\tW$ is attained for those $r$. 

We now fix such an $r > 1$ near $1$ and choose $F = \Id + G$ such that $\tW(r) = W_{\triangle}(F)$ and $\e_1^T F \e_1 = r$. As $W_{\triangle}$ is invariant under the reflection $P$, we may without loss of generality assume that $G$ is small. Then Lemma \ref{lemma:W-linearized-properties} yields 
$$ W_{\triangle}(F) 
   = \frac{3 \alpha}{32} \left( 3 g_{11}^2 + 3 g_{22}^2 + 2 g_{11} g_{22} 
     + 4 \left(\frac{g_{12}+g_{21}}{2} \right)^2 \right) + o(|G|^2). $$
We find that $g_{11} = r - 1$, $g_{12} + g_{21} = o(r-1)$ and $g_{22} = - \frac{1}{3} g_{11} + o(r-1)$ and $F$ satisfies 
$$ \frac{F^T + F}{2} = \begin{pmatrix} r & 0 \\ 0 & \frac{4 - r}{3} \end{pmatrix} + o(r-1) $$ 
with energy 
\begin{align*}
  W_{\triangle} (F) 
  &= W_{\triangle}\left( \frac{F^T + F}{2} \right) + o((r-1)^2) \\ 
  &= \frac{\alpha}{4} \left( r - 1 \right)^2 + o((r-1)^2). 
\end{align*} 

(iii) This is immediate from Lemma \ref{lemma:W-triangle-properties}(iii). \eop

Under strengthened hypotheses on $W$ we have the following expansion:
\begin{lemma}\label{lemma:Wtilde-fein} 
If $W$ in addition satisfies the assumptions (ii') and (iii'), then for $r > 1$ close to $1$ we have 
$$ \tilde{W}(r) 
   = \frac{\alpha (r-1)^2}{4} 
     + \frac{1}{108} \Big( 6 \alpha +  7 \alpha'  - 2 ( 3 \alpha - \alpha' ) \cos (6 \phi) \Big) (r-1)^3 + O((r-1)^4), $$
where $\phi$ is such that $R_{\cal L} = \begin{footnotesize} \begin{pmatrix} \cos \phi & - \sin \phi \\ \sin \phi & \cos \phi \end{pmatrix} \end{footnotesize}$
\end{lemma} 

\Proof Let $s = r-1$. By definition, 
$$ \tilde{W}(r) 
   = \min \left\{ W_{\triangle} (F(s, x, y, z)) : x, y, z \in \R \right\}, $$
where $F(s, x, y, z) =  \begin{footnotesize} \begin{pmatrix} 1 + s & z + y \\ z - y & 1 + x \end{pmatrix} \end{footnotesize}$. Due to the quadratic energy growth near $SO(2)$, we need to minimize only over $x, y, z$ with $|x|, |z|, \sqrt{s} |y| \le C s$ for a constant $C$ large enough. Indeed, as $W_{\triangle} (F(s, 0, 0, 0)) = O(s^2)$, for a minimizer one has without loss of generality $\dist(F(s, x, y, z), SO(2)) = O(s)$. But then $\sqrt{(1 + s)^2 + (z \pm y)^2} = 1 + O(s)$, which implies $|z \pm y| = O(\sqrt{s})$ and so $|z|, |y| = O(\sqrt{s})$, and also $\sqrt{(1 + x)^2 + (z \pm y)^2} = 1 + O(s)$, which then implies $\pm (1 + x) = 1 + O(s)$ and thus without loss of generality $x = O(s)$. Finally using that the scalar product $(1 + s)(z + y) + (1 + x)(z - y) = 2z + O(s^{3/2})$ of the two columns of $F(s, x, y, z)$ in absolute value is also bounded by $O(s)$, we obtain that $|z| = O(s)$. 

Set $x = - \frac{s}{3} + s x_1$, $y = \sqrt{s} y_1$, $z = s z_1$ with $|x_1|, |y_1|, |z_1| \le C$. Explicit calculation gives 
\begin{align*}
  W_{\triangle}(F(s, x, y, z)) 
  = \frac{\alpha}{32} \left( 8 + 3 x_1^2 + 8 y_1^2 + 12 z_1^2 + 6 ( x_1 + y_1^2 )^2 \right) s^2 
     + O(s^3). 
\end{align*}
Since $\alpha > 0$, we thus obtain that this expression is minimized in $x_1, y_1, z_1$ with $x_1^2, y_1^2, z_1^2 = O(s)$ and we may set $x_1 = \sqrt{s} x_2$, $y_1 = \sqrt{s} y_2$ and $z_1 = \sqrt{s} z_2$ with $|x_2|, |y_2|, |z_2| \le C$ for some $C > 0$. Explicit expansion in powers of $s$ then yields 
\begin{align*}
  &W_{\triangle}(F(s, x, y, z)) \\ 
  &= \frac{\alpha s^2}{4} 
     + \frac{1}{864} \Big( 48 \alpha + 56 \alpha' - 16 (3 \alpha - \alpha') \cos (6 \phi) \\ 
     &\qquad \qquad \qquad \qquad \qquad 
     + 3 \alpha \left( 81 x_2^2 + 72 y_2^2 + 108 z_2^2 \right) \Big) s^3 \\ 
  &\quad   + \frac{1}{24} \Big( \left(9 \alpha y^2 + \alpha' + (3 \alpha - \alpha') \cos( 6 \phi ) \right) x_2 \\ 
     &\qquad \qquad \qquad \qquad \qquad 
    + 2( 3 \alpha - \alpha') \sin( 6 \phi ) z_2) \Big) s^{7/2} 
     + O(s^4) \\ 
  &= \frac{\alpha s^2}{4} 
     + \frac{1}{108} \Big( 6 \alpha +  7 \alpha'  - 2 ( 3 \alpha - \alpha' ) \cos (6 \phi) \Big)s^3 \\ 
     &\quad +  \frac{9 \alpha}{32} \Big( x_2^2 + 2 A \sqrt{s} x_2 \Big) s^3 
     + \frac{\alpha y_2^2 s^3}{4} 
     + \frac{3 \alpha}{8} \Big(z_2^2 + 2 B \sqrt{s} z_2 \Big) s^3  + O(s^4) 
\end{align*}
for $A$ and $B$ bounded uniformly in $s$ and so 
\begin{align*}
  &W_{\triangle}(F(s, x, y, z)) \\ 
  &= \frac{\alpha s^2}{4} 
     + \frac{1}{108} \Big( 6 \alpha +  7 \alpha'  - 2 ( 3 \alpha - \alpha' ) \cos (6 \phi) \Big)s^3 \\ 
     &\quad +  \frac{9 \alpha}{32} \Big( x_2 + A \sqrt{s} \Big)^2 s^3 
     + \frac{\alpha y_2^2 s^3}{4} 
     + \frac{3 \alpha}{8} \Big(z_2 + B \sqrt{s} \Big)^2 s^3  + O(s^4). 
\end{align*}
Minimizing with respect to $x_2, y_2$ and $z_2$ we finally obtain that 
$$ \tilde{W}(1+s) 
   = \frac{\alpha s^2}{4} 
     + \frac{1}{108} \Big( 6 \alpha +  7 \alpha'  - 2 ( 3 \alpha - \alpha' ) \cos (6 \phi) \Big)s^3 + O(s^4). $$
\eop 

The following lemma provides useful lower bounds for the energy $W_{\triangle}$ and the reduced energy $\tW$. 
\begin{lemma}\label{lemma:quadratic lower bound} 
For all $T>1$ one has:
\begin{itemize}
\item[(i)] There exists some $c > 0$ such that $c \dist^2(F,O(2)) \leq  W_{\triangle} (F)$ for all $F \in \R^{2 \times 2}$ satisfying $|F|\leq T$.
\item[(ii)] For $\delta > 0$ small enough, there is a convex function $V \geq 0$ with $V(r) \le \tilde{W}(r)$ for $r \le T$ and such that the second derivative $V''_+(1)$ from the right at $1$ exists and satisfies $V''_+(1) = \frac{\alpha}{2} - 2\delta$. 
\item[(iii)] If in addition $W$ satisfies assumptions (ii') and (iii'), then there exists a convex function $V \ge 0$ with $V(r) \le \tilde{W}(r) \le V(r) + O((r-1)^4)$ for $r \le T$. 
\end{itemize}
\end{lemma}

\Proof
(i) Let $F \in \R^{2 \times 2}$ satisfying $|F| \leq T$. By polar decomposition we find $R \in O(2)$ and $U=\sqrt{F^{T}F}$ symmetric and positiv definite such that $F = RU$. A short computation yields $|U - \Id| = \dist(F,O(2))$. Assume first $|U - \Id| < \eta$ for $\eta > 0$ small enough. Since $W_{\triangle}(F) $ is invariant under rotation and reflection we obtain applying Lemma \ref{lemma:W-linearized-properties}: 
$$W_{\triangle}(F) = W_{\triangle}(R^{T}RU)\geq \frac{1}{2}Q(U - \Id) + o(|U - \Id|^2).$$
Noting that $Q$ grows quadratically on $\R^{2 \times 2}_{\text{sym}}$ (see Lemma \ref{lemma:W-linearized-properties}) we obtain a constant $c_{1} > 0$ such that for $|U - \Id| < \eta$ 
$$ W_{\triangle}(F) \geq c_1 |U - \Id|^2 = c_1 \dist^2(F,O(2)).$$

Consider the compact set $M:=\left\{F \in \R^{2 \times 2}, \dist(F,O(2))\geq\eta, |F|\leq T\right\}$. $W_{\triangle}$ attains its minimum on M, which is strictly positiv by Lemma \ref{lemma:W-triangle-properties}(ii). This provides a second constant $c_2 > 0$  such that for all $F \in M$
$$ W_{\triangle}(F) \geq c_2 |U - \Id|^2 = c_2 \dist^2(F,O(2)).$$ 
Taking $c=\min\{c_1,c_2\}$ yields the claim. 

(ii) We construct such a function directly applying Lemma \ref{lemma:W-reduced-properties}.
$$ V(r) 
= \begin{cases} 
        0 
          & \text{for } r \le 1, \\ 
        \left(\frac{\alpha}{4} - \delta\right)(r - 1)^2 
          & \text{for } 1 \le r \le 1 + \eta, \\ 
        \left(\frac{\alpha}{4} - \delta\right)\eta \left(2r -2 -\eta\right)
          & \text{for } r \ge 1 + \eta, 
     \end{cases} $$ 
when $\eta > 0$ is sufficiently small. 

(iii) With $f(r) := \frac{\alpha (r-1)^2}{4} + \frac{1}{108} \Big( 6 \alpha +  7 \alpha'  - 2 ( 3 \alpha - \alpha' ) \cos (6 \phi) \Big) (r-1)^3 - C(r-1)^4$ for sufficiently large $C$, Lemma \ref{lemma:Wtilde-fein} shows that we can choose 
$$ V(r) 
= \begin{cases} 
        0 
          & \text{for } r \le 1, \\ 
        f(r) 
          & \text{for } 1 \le r \le 1 + \eta, \\ 
        f(1+\eta) + f'(1+\eta) (r - 1 - \eta) 
          & \text{for } r \ge 1 + \eta, 
     \end{cases} $$ 
when $\eta > 0$ is sufficiently small. 
\eop

\section{Limiting minimal energy and cleavage laws}\label{sec:limiting-energy}

We now prove Theorems \ref{theo:limiting-energy} and \ref{theo: discrete energy} on cleavage laws and fine energy estimates. 

\subsection*{Limiting minimal energy}

We can classify (or `color') all triangles in ${\cal C}_{\eps}$ into two types, say `type one' and `type two', such that all triangles of the same type are translates of each other. Then only triangles of different type can share a common side. Denote the sets by ${\cal C}^{(1)}_{\eps}$ and ${\cal C}^{(2)}_{\eps}$, respectively. \smallskip

\noindent {\em Proof of Theorem \ref{theo:limiting-energy}.}  We first show that the expression on the right hand side is a lower bound for the limiting minimal energy. For every deformation $y \in {\cal A}(a_{\eps})$ we have by \eqref{eq:E-integral} and \eqref{eq:W-triangle} 
\begin{align*}
  {\cal E}_{\eps} (y) 
  &\ge \frac{4}{\sqrt{3}\eps} \int_{\Omega_{\eps} \cap (0,l) \times (\eps, 1 - \eps)} W_{\triangle} \big( \nabla \tilde{y}\big) \, dx. 
\end{align*}

Let $0 < \delta < \frac{\alpha}{4}$ and choose $R$ so large that $W(r) > \beta - \delta$ if $r \ge R$. Define $\bar{\cal C}^{(1)}_{\eps}$ to be the set of those triangles $\triangle$ of type one for which at least one side in the deformed configuration $y(\triangle)$ is larger than $2R\eps$. By $I \subset (\eps, 1 - \eps)$ we denote the set of those points $x_2$ for which there exists $x_1 \in (0, l)$ such that $(x_1, x_2)$ lies in one of these triangles. 

We can then estimate the energy integral by splitting the $x_2$-integration into a first part where $x_2 \notin I$ and a second part with $x_2 \in I$. 

1. If $x_2 \notin I$, then all sidelengths of $y(\triangle)$ for a triangle $\triangle$ whose interior intersects the segment $(0, l) \times \{x_2\}$ are less or equal to $4R\eps$. This is clear for triangles of type one by construction. For triangles of type two it follows from the fact that the two sides of $\triangle$ intersecting $(0,l) \times \{x_2\}$ are also sides of triangles of type one and therefore bounded by $2R\eps$. The third side is thus less than $4R\eps$, too. 

It is elementary to see that for $F \in \R^{2 \times 2}$  
\begin{equation}\label{eq: bondlengths}
|\e_1^{T} F \e_1| \le 8R \text{,  \ \ \ if } |\vv^T F \vv| \le 4R \text{ for all } \vv \in {\cal V}. 
\end{equation}
Indeed, if $\lambda_1, \lambda_2$ are the eigenvalues of $\frac{1}{2}(F^T + F)$, then by \eqref{eq:lat-sum} one has $\frac{3}{4}(\lambda_1^2 + \lambda_2^2) = \frac{3}{4} \trace\left(\frac{1}{2}(F^{T} + F)\right)^2 \le 3 \cdot (4R)^2$ and thus $|\e_1^{T} F \e_1| \le \max\{|\lambda_1|, |\lambda_2|\} \le 8R$. Consequently, for almost every $x_2 \notin I$ we have $\e_1^T \nabla \tilde{y}(x_1, x_2) \e_1 \le 8R$ for all $x_1 \in (0, l)$. 

By Lemma \ref{lemma:quadratic lower bound}(ii) choose a convex function with $V(r) \le \tilde{W}(r)$ for $r \le 8R$ and $V''_+(1) = \frac{ \alpha}{2} - 2\delta$. For $x_2 \in (\eps, 1 - \eps)$ define $\Omega^{x_2}_{\eps} \subset (0,l)$ such that $\Omega^{x_2}_{\eps} \times \{ x_2\} = \Omega_{\eps} \cap (0,l) \times \{ x_2\}$. Then for the first part one obtains, if $a < \infty$, by convexity of $V$
\begin{align}\label{eq: elastic energy estimate} 
\begin{split}
  \frac{4}{\sqrt{3}\eps} \int_{(\eps,1 -\eps)\setminus I} \int_{\Omega^{x_2}_{\eps}}W_{\triangle} \big( \nabla \tilde{y}\big) \, dx_1 \, dx_2 
  &\ge \frac{4}{\sqrt{3}\eps} \int_{(\eps,1-\eps)\setminus I} 
   \int_{\Omega^{x_2}_{\eps}} V \big( \e_1^T \nabla \tilde{y} \, \e_1 \big) \, dx_1 \, dx_2 \\ 
  &\ge \frac{4}{\sqrt{3}\eps} \int_{(\eps,1-\eps)\setminus I} |\Omega^{x_2}_{\eps}| V(1 + a_{\eps}) \, dx_2 \\ 
  &= \frac{2}{\sqrt{3}\eps} (1 - 2\eps - |I|) (l- 2\eps) ( V''_+(1) a_{\eps}^2 + o(\eps) ) \\ 
  &\to \frac{2}{\sqrt{3}} (1 - |I|) l V''_+(1) a^2 
\end{split}
\end{align}
as $\eps \to 0$. It is not hard to see that this asymptotic estimate remains true also for $a = \infty$. 

2. On the other hand, the energy of the second part can be estimated by the energy of all springs lying on the side of a triangle in $\bar{\cal C}^{(1)}_{\eps}$, which yields 
\begin{align}\label{eq: crack energy estimate}
  \frac{4}{\sqrt{3}\eps} \int_{I} \int_{\Omega^{x_2}_{\eps}} W_{\triangle} \big( \nabla \tilde{y}\big) \, dx_1 \, dx_2 
  \ge 2 (\beta - \delta) \eps \#\bar{\cal C}^{(1)}_{\eps}, 
\end{align}
as the length of at least two springs in each of these triangles is larger than $R\eps$ in the deformed configuration. Now the projection of any triangle onto the $x_2$-axis is an interval of length $\eps \gamma$, and so $\eps \gamma \#\bar{\cal C}^{(1)}_{\eps} \ge |I|$, i.e., 
\begin{align}\label{eq:crack-energy-interval}
  \frac{4}{\sqrt{3}\eps} \int_{I} \int_{\Omega^{x_2}_{\eps}} W_{\triangle} \big( \nabla \tilde{y}\big) \, dx_1 \, dx_2 
  \ge 2 (\beta - \delta) \gamma^{-1} |I|. 
\end{align}

Summarizing \eqref{eq: elastic energy estimate} and \eqref{eq:crack-energy-interval} we find 
\begin{align*}
  &\liminf_{\eps \to \infty} \min \{ {\cal E}_{\eps}(y) : y \in {\cal A}(a_{\eps}) \} \\ 
  &\ge \min \left\{\frac{2}{\sqrt{3}} \left(\frac{\alpha}{2} - 2 \delta\right)  l a^2 (1 - |I|) 
       + 2 (\beta - \delta) \gamma^{-1} |I| : |I| \in [0, 1] \right\} \\ 
  &= \min \left\{\frac{2}{\sqrt{3}} \left(\frac{\alpha}{2} - 2 \delta\right) l  a^2, 
       \frac{2(\beta - \delta)}{\gamma} \right\}. 
\end{align*}Now $\delta \to 0$ shows 
\begin{align*}
  \liminf_{\eps \to \infty} \min \{ {\cal E}_{\eps}(y) : y \in {\cal A}(a_{\eps}) \} 
  \ge \min \left\{\frac{\alpha l}{\sqrt{3} } a^2, \frac{2 \beta}{\gamma} \right\}. 
\end{align*}
This establishes the lower bound. 

It remains to prove that the right hand side in Theorem \ref{theo:limiting-energy} is attained for some sequence of deformations. In order to do so, we consider two specific sequences of deformations. First, for $a < \infty$ let 
\begin{align}\label{eq: elastic minimizer}
  y^{\rm el}_{\eps}(x) 
  = (\Id + F^{a_{\eps}}) x 
  = \begin{pmatrix} 
        1 + a_{\eps} & 0 \\ 
        0 & 1 - \frac{a_{\eps}}{3} 
      \end{pmatrix} x. 
\end{align}
By Lemma \ref{lemma:W-reduced-properties}(ii) we have that $W_{\triangle}(F) = \frac{\alpha}{4} a_{\eps}^2 + o(\eps)$ and so 
\begin{align*}
  \lim_{\eps \to 0} {\cal E}_{\eps}(y^{\rm el}_{\eps}) 
  = \frac{\alpha l}{\sqrt{3} } a^2 
\end{align*}
by \eqref{eq:E-integral}. 

To define $y^{\rm cr}$ we choose any line $(s,0) + \R\vv_{\gamma}$ intersecting both the segments $(0, l) \times \{0\}$ and $(0, l) \times \{1\}$ (as in Corollary \ref{cor: limiting-deformations}). This is possible since $l > \frac{1}{\sqrt{3}}$.  Let $a > 0$ and set 

\begin{align}\label{eq: crack minimizer}
  y^{\rm cr}_{\eps}(x) 
  = \begin{cases} 
      x & \mbox{for } x \mbox{ to the left of } (s,0) + \R\vv_{\gamma}, \\ 
      x + a_{\eps} l \e_1 & \mbox{for } x \mbox{ to the right of } (s,0) + \R\vv_{\gamma} 
    \end{cases} 
\end{align} 
for atoms $x$ with $\eps < x_1 < l-\eps$. Except for a negligible contribution from the boundary layers, the energy of this configuration can be estimated as in Step 2 of the proof of the lower bound: It is given by the energy of springs intersecting $(s,0) + \R\vv_{\gamma}$, i.e., by the two springs lying on the boundary of the triangles of type one which are intersected by $(s,0) + \R\vv_{\gamma}$. These springs are elongated by a factor scaling with $a_{\eps}/\eps$, thus yielding a contribution $\beta$ in the limit $\eps \to 0$.  \eop

\subsection*{Fine estimates on the limiting minimal energy}

Assume now that $W$ in addition satisfies assumptions (ii') and (iii'). In order to investigate a deformation $y$ again we let $\bar{{\cal C}}_{\eps}$ and $\bar{{\cal C}}_{\eps}^{(1)}$ denote the set of triangles $\triangle$ (of type one respectively) for which at least one side in $y(\triangle)$ is larger than $2R\eps$, where now the threshold value $R > 1$ is chosen in such a way that $c_R := \inf\{W(r) : r \ge R\} \ge \frac{\beta}{2}$. According to Lemma \ref{lemma:quadratic lower bound}(iii) we may choose a convex function $V$ such that 
\begin{align}\label{eq:V-choice} 
  0 \le V(r) \le \tilde{W}(r) \le V(r) + O((r-1)^4) \mbox{ for } r \le 8R. 
\end{align} 
As in \eqref{eq: bondlengths} we observe that $|\e_1^{T} (y)_{\triangle} \e_1|$ is bounded by $8R$ on  triangles with bond length not exceeding $4R\eps$ and thus lies in the convex regime of $V$. Moreover, we find that every triangle in $\bar{{\cal C}}_{\eps}$ provides at least the energy $\frac{4}{\sqrt{3}\eps} \int_{\triangle} W_{\triangle}(\nabla \tilde{y}) \ge c_R \eps$. 

For given $0 < \eta < a$ we also define $R_{\eps, \eta} = \frac{a-\eta}{\sqrt{\eps}}$ as a threshold for triangles we consider `essentially broken': 
\begin{equation}\label{eq: real, pseudo}
  \bar{{\cal C}}_{\eps,\eta} 
  = \left\{\triangle \in \bar{{\cal C}_{\eps}}, |\nabla y_{\eps} \vv| > R_{\eps,\eta} 
       \text{ for at least two } \vv \in {\cal V}\right\}. 
\end{equation}
The minimal energy contribution of all the springs on such a triangle in $\bar{{\cal C}}_{\eps,\eta}$ is given by 
$$ 2 \beta^{\eta} \eps
   := 2 \inf \left\{W(r) : r \ge \frac{a-\eta}{\sqrt{\eps}} \right\} \eps
   = (2 \beta + O(\eps)) \eps $$ 
by the assumption (iii') on $W$. By $I \subset (\eps, 1-\eps)$ we denote the set of points $x_2$ for which the segment $(0,l) \times \left\{x_2\right\}$ intersects a broken triangle (of type one) in $\bar{{\cal C}}^{(1)}_{\eps}$. In addition, we say $x_2 \in I^{\eta} \subset I$ if one of the intersected triangles lies in $\bar{{\cal C}}_{\eps,\eta} \cap \bar{{\cal C}}^{(1)}_{\eps}$. 

With these preparations we can now proceed to prove Theorem \ref{theo: discrete energy}: \smallskip 

\noindent{\em Proof of Theorem \ref{theo: discrete energy}.} 
Let ${\cal E}_{\eps}(y) = \inf{\cal E}_{\eps} + O(\eps)$. Inspired by (\ref{eq: elastic energy estimate}) and (\ref{eq: crack energy estimate}) we establish a lower bound for the energies additionally taking the set $I \setminus I^{\eta}$ into account. Since the sidelength of any triangle whose interior intersects $(0, l) \times (I \setminus I^{\eta})$ is bounded by $4R_{\eps,\eta}$, we find 
$$ | \e_1^{T} \nabla \tilde{y}(x_1, x_2) \, \e_1 | 
   \le 8 R_{\eps, \eta} $$
for all $(x_1,x_2) \in (0, l) \times (I \setminus I^{\eta})$ as in \eqref{eq: bondlengths}. Let $k = k(x_2)$ count the number of triangles in $\bar{\cal C}_{\eps}$ on the slice $(0, l) \times \{x_2\}$, $x_2 \in I \setminus I^{\eta}$, and define $\bar{\cal C}^{x_2}_{\eps} \subset (0,l)$ such that $((0, l) \times \{x_2\}) \cap \bigcup_{\triangle \in \bar{{\cal C}}_{\eps}} \triangle = \bar{\cal C}^{x_2}_{\eps} \times \{x_2\}$. Then 
$$ \int_{\bar{\cal C}^{x_2}_{\eps}} \e_1^{T} \nabla \tilde{y}(x_1, x_2) \, \e_1 \, dx_1
   \le 8 k \eps  R_{\eps, \eta}. $$
and so 
\begin{align*} 
  \int_{\Omega^{x_2}_{\eps} \setminus \bar{\cal C}^{x_2}_{\eps}} \e_1^{T} \nabla \tilde{y}(x_1, x_2) \, \e_1
  &\ge (1 + \sqrt{\eps} a) (l + O(\eps)) - 8 k \eps  R_{\eps, \eta} \\ 
  &= \left( 1 + \sqrt{\eps} \left( a - \frac{8 k (a - \eta)}{l} + O(\sqrt{\eps}) \right) \right) l . 
\end{align*}
Since $\#(\bar{{\cal C}}_{\eps} \setminus \bar{{\cal C}}_{\eps,\eta}) \ge \frac{1}{\eps \gamma} \int_{I \setminus I^{\eta}} k(x_2) \, dx_2$, a convexity argument as in the proof of Theorem \ref{theo:limiting-energy} on slices $(0, l) \times \{x_2\}$ with $x_2 \in (\eps,1 - \eps) \setminus I$ and on the unbroken part $\big(\Omega^{x_2}_{\eps} \setminus \bar{\cal C}^{x_2}_{\eps}\big) \times \{x_2\}$ of slices with $x_2$ in $I \setminus I^{\eta}$ then shows that 
\begin{align} \label{eq: lower bound, three slice types}
  {\cal E}_{\eps}(y) 
  \geq \frac{4 (l - 2\eps)}{\sqrt{3} \eps} V(1 + \sqrt{\eps}a) (1 - 2\eps - |I|) 
       +  G_{\eta,\eps} |I \setminus I^{\eta}|
       + \frac{2 \beta^{\eta}}{\gamma} |I^{\eta}| + O(\eps),  
\end{align}
where
\begin{align*}
  G_{\eta,\eps}
  = \min_{k \in \N} \left(\frac{4 l}{\sqrt{3} \eps} 
    V\left(1 + \sqrt{\eps}\left(a - \frac{8k(a - \eta)}{l} + O(\sqrt{\eps})\right)\right) 
    + \frac{k c_R}{\gamma} \right).  
\end{align*}
We note that this minimum exists and can be taken over $1 \leq k \leq K_{0}$ for some $K_{0} \in \N$ large enough and independent of $\eta$ as $\frac{k c_R}{\gamma} \to \infty$ for $k \rightarrow \infty$. We choose $0 < \eta < a$ large enough such that 
$$ \frac{ l \alpha}{\sqrt{3}} a^2 
   < \min_{1 \leq k \leq K_0} 
     \left( \frac{ \alpha l}{\sqrt{3}} \left(a - \frac{8k}{l} (a-\eta) \right)^2 + \frac{k c_R}{\gamma} \right). $$ 
Recalling that, by (\ref{eq:V-choice}) and Lemma \ref{lemma:W-reduced-properties}, $\frac{4 l}{\sqrt{3} \eps}V(1 + \sqrt{\eps} r) = \frac{4 l}{\sqrt{3} \eps}\tW(1 + \sqrt{\eps} r) + O(\eps) \rightarrow \frac{ l \alpha}{\sqrt{3}} r^2$ uniformly in $r$ on bounded sets in $\R$, we see that thus $G_{\eta,\eps}$ exceeds the elastic term $\frac{4 l}{\sqrt{3} \eps}V(1 + \sqrt{\eps}a)$ for $\eps$ sufficiently small. So from \eqref{eq: lower bound, three slice types} we obtain
\begin{align}\label{eq: lower bound, two slice types}
  {\cal E}_{\eps}(y) 
  \geq \frac{4 l}{\sqrt{3} \eps} V(1 + \sqrt{\eps}a) (1 - 2\eps - |I^{\eta}|) 
   + \frac{2 \beta^{\eta}}{\gamma} |I^{\eta}| + O(\eps). 
\end{align}
As $\frac{4 l}{\sqrt{3} \eps}V(1 + \sqrt{\eps}a) \to \frac{l \alpha}{\sqrt{3}} a^2$ and $\beta^{\eta} \rightarrow \beta$ for all $\eta>0$, for $\eps$ small enough we thus obtain $\inf{\cal E}_{\eps} \geq \frac{4}{\sqrt{3}}\frac{l}{\eps} V(1+\sqrt{\eps}a) (1-2\eps) + O(\eps) = \frac{4}{\sqrt{3}}\frac{l}{\eps} \tW(1+\sqrt{\eps}a) + O(\eps)$ or $\inf{\cal E}_{\eps} \geq \frac{2 \beta + O(\eps)}{\gamma}(1-2\eps) = \frac{2 \beta}{\gamma} + O(\eps)$, respectively, depending on $a$. 

Applying (\ref{eq: elastic minimizer}) and (\ref{eq: crack minimizer}) we then get indeed $\inf{\cal E}_{\eps} = \frac{4}{\sqrt{3}}\frac{l}{\eps} \tW(1+\sqrt{\eps}a) + O(\eps)$ or $\inf{\cal E}_{\eps} = \frac{2 \beta}{\gamma} + O(\eps)$, respectively. The claim now follows from Lemma \ref{lemma:Wtilde-fein} . 
\eop

\noindent{\bf Remark.} From the proof of Theorem \ref{theo: discrete energy}, especially taking (\ref{eq: elastic minimizer}) and (\ref{eq: crack minimizer}) into account, it follows that Theorem \ref{theo: discrete energy} still holds if ${\cal E}_{\eps}$ is replaced by ${\cal E}^{\chi}_{\eps}$.

\section{Limiting minimal energy configurations}\label{sec:limiting-configurations}

Throughout this section we will assume that $a_{\eps} = \sqrt{\eps} a$, $y_{\eps}$ is a sequence of deformations satisfying \eqref{eq:almost-min}, the threshold value $R$ is chosen as above Equation \eqref{eq:V-choice} and that $\bar{\cal C}_{\eps}$ is defined accordingly. 

For a rescaled displacement $\tilde{u}$ we denote by $D^{\mu} \subset (\eps,1-\eps)$ for $\mu > 0$ the set of $x_2$ such that there is precisely one triangle $\triangle_{x_2} \in \bar{\cal C}^{(1)}_{\eps}$ with $\text{int}(\triangle_{x_2}) \cap \left((0, l) \times \{x_2\}\right) \neq \emptyset$ and 
\begin{equation}\label{eq: D_mu}
\int_{\Omega^{x_2}_{\eps} \setminus \bar{\cal C}^{x_2}_{\eps}} \e_1^{T}\nabla \tilde{u}(x_1, x_2) \e_1 \, dx_1  \leq l\mu.
\end{equation}
Note that $D^{\mu} \subset I^{\eta}$ for $\mu$ small enough: For $x_2 \in D^{\mu}$ we have 
$$\int_{\bar{\cal C}^{x_2}_{\eps}} \e_1^{T}\nabla \tilde{y}(x_1, x_2) \e_1 \, dx_1 \geq \sqrt{\eps} l (a - \mu)  + O(\eps)$$
and using the arguments in (\ref{eq: bondlengths}) we see that for given $\eta$ (not too small) we can choose $\mu$ small enough such that $\triangle_{x_2} \in \bar{\cal C}_{\eps,\eta}$ and thus $x_2 \in I^{\eta}$. We also define $\bar{{\cal C}}^{\mu}_{\eps,\eta} \subset \bar{{\cal C}}_{\eps,\eta}$ as the set of those essentially broken triangles $\triangle$ for which there exists some $x_2 \in D^{\mu}$ such that $\text{int}\left(\triangle\right) \cap \left((0,l) \times \left\{x_2\right\}\right) \neq \emptyset$. The projection of a triangle $\triangle$ onto the linear subspace spanned by $\vv^{\bot}_{\gamma}$ is an interval of length $\frac{\sqrt{3}}{2}\eps$. We denote the center of this interval by $m_{\triangle}$.

The following lemmas give sharp estimates on the number of broken triangles and their position. 

\begin{lemma}\label{lemma: minimizing sequences - properties - elast}
Let $a < a_{\rm crit}$ and suppose $\tilde{u}_{\eps}$ is a minimizing sequence satisfying 
$${\cal E}_{\eps}(\id + \sqrt{\eps} u_{\eps}) = \inf {\cal E}_{\eps}+O(\eps).$$ 
Then $\eps \# \bar{{\cal C}}_{\eps} = O(\eps)$. 
\end{lemma} 

\Proof
Using \eqref{eq: lower bound, three slice types} we find 
\begin{align*}
  {\cal E}_{\eps}(y_{\eps}) 
  &= \frac{4 l}{\sqrt{3} \eps} \tW(1 +\sqrt{\eps}a) + O(\eps) \\ 
  &\geq \frac{4 (l-2\eps)}{\sqrt{3} \eps}\tW(1 + \sqrt{\eps}a) (1 - 2\eps - |I|) 
   + \min\left\{G_{\eta,\eps}, \frac{2 \beta^{\eta}}{\gamma}\right\} |I| + O(\eps) \\
   & = \frac{4 l}{\sqrt{3} \eps}\tW(1 + \sqrt{\eps}a) (1 -  |I|) 
   + \min\left\{G_{\eta,\eps}, \frac{2 \beta^{\eta}}{\gamma}\right\} |I| + O(\eps).
\end{align*}
An elementary computation yields, whenever $\eps$ is small enough,
\begin{align*}
  |I|  
  &\leq \left(\min\left\{G_{\eta,\eps}, \frac{2 \beta^{\eta}}{\gamma}\right\} 
    - \frac{4 l}{\sqrt{3} \eps} \tW(1 +\sqrt{\eps}a)\right)^{-1} \cdot O(\eps) \\
  &= \left(\min\left\{G_{\eta,\eps}, \frac{2 \beta}{\gamma}\right\} - \frac{\alpha l}{\sqrt{3} }a^2 + o(1)\right)^{-1} \cdot O(\eps)             = O(\eps). 
 \end{align*}
(The argument leading to \eqref{eq: lower bound, two slice types} together with $a < a_{\rm crit}$ shows that the term in parentheses is bounded from below by a positive constant independent of $\eps$). Then the elastic energy is $\frac{4 l}{\sqrt{3} \eps} \tW(1 +\sqrt{\eps}a) + O(\eps)$ and consequently, the crack energy coming from triangles in $\bar{{\cal C}}_{\eps}$ is of order $O(\eps)$. As every broken triangle in $\bar{{\cal C}}_{\eps}$ provides at least energy $\eps c_R$ we conclude $\eps \# \bar{{\cal C}}_{\eps}^{(1)} = O(\eps)$. But then, possibly after replacing $R$ by $2R$, also $\eps \# \bar{{\cal C}}_{\eps}^{(2)} = O(\eps)$ as those triangles are neighbors of broken triangles of type $1$.
\eop 

\begin{lemma}\label{lemma: minimizing sequences - properties - crack}
Let $a > a_{\rm crit}$, $\phi \ne 0$ and suppose $\tilde{u}_{\eps}$ is a minimizing sequence satisfying 
$${\cal E}_{\eps}(\id + \sqrt{\eps} u_{\eps}) = \inf {\cal E}_{\eps}+O(\eps).$$ 
Then $|I^{\eta}| = 1- O(\eps)$ for $0 < \eta < a$. Furthermore, for $\mu$ sufficiently small, $\eps \# \left(\bar{{\cal C}}_{\eps} \setminus \bar{{\cal C}}^{\mu}_{\eps,\eta}\right) = O(\eps)$ and 
$$\sup \left\{  |m_{\triangle_{1}} - m_{\triangle_{2}}| : \triangle_{1}, \triangle_{2} \in \bar{{\cal C}}^{\mu}_{\eps,\eta} \right\} = O(\eps).$$
\end{lemma} 

\Proof Using (\ref{eq: lower bound, two slice types}) we find after without loss of generality choosing $\eta$ sufficiently large 
\begin{align*}
  {\cal E}_{\eps}(y_{\eps}) 
  = \frac{2 \beta}{\gamma} + O(\eps) 
  \geq \frac{4 l}{\sqrt{3} \eps} \tW(1 + \sqrt{\eps}a) (1 - 2\eps - |I^{\eta}|) 
   + \frac{2 \beta^{\eta}}{\gamma} |I^{\eta}|.
\end{align*}
So for $\eps$ small enough we obtain
$$ 1- |I^{\eta}| 
   \leq \left(\frac{ \alpha l}{\sqrt{3}} a^2 + o(1) - \frac{2 \beta}{\gamma}\right)^{-1} \cdot O(\eps) 
   = O(\eps) $$
since $a > a_{\rm crit}$. Consequently, the crack energy from triangles in $\bar{{\cal C}}_{\eps,\eta}$ is given by $\frac{2\beta}{\gamma} + O(\eps)$ and thus the energy contribution from $\bar{{\cal C}}_{\eps} \setminus \bar{{\cal C}}_{\eps,\eta}$ is of order $O(\eps)$. As in the proof of Lemma \ref{lemma: minimizing sequences - properties - elast} we find $\eps \# \left(\bar{{\cal C}}_{\eps} \setminus \bar{{\cal C}}_{\eps,\eta}\right) = O(\eps)$. 
Let $k_{\eta}(x_2)$ and $k^C_\eta(x_2)$ count the number of triangles in $\bar{{\cal C}}_{\eps,\eta} \cap \bar{{\cal C}}_{\eps}^{(1)}$ and $(\bar{{\cal C}}_{\eps} \setminus \bar{{\cal C}}_{\eps,\eta}) \cap \bar{{\cal C}}_{\eps}^{(1)}$ intersected by $(0,l) \times \{x_2\}$, respectively. We dissect $I^{\eta} \setminus D^{\mu}$ into two disjoint sets: By $D_{1} \subset I^{\eta} \setminus D^{\mu}$ we denote the set where we find more than one triangle $\triangle_{x_2} \in \bar{\cal C}^{(1)}_{\eps}$ with $\text{int}(\triangle_{x_2}) \cap \left((0, l) \times \{x_2\}\right) \neq \emptyset$. The complement $D_2$ is the set where \eqref{eq: D_mu} does not hold. Using a convexity argument for $x_2 \in D_2$ we obtain 
\begin{align*}
  \frac{2\beta}{\gamma} + O(\eps) 
  & \geq 2 \beta \eps (\#\bar{{\cal C}}_{\eps,\eta} \cap \bar{{\cal C}}_{\eps}^{(1)}) 
     + 2 c_R \eps \# \left( (\bar{{\cal C}}_{\eps} \setminus \bar{{\cal C}}_{\eps,\eta}) \cap \bar{{\cal C}}_{\eps}^{(1)} \right) \\  
  &\qquad + \frac{4 }{\sqrt{3} \eps} \int^{1-\eps}_\eps \int_{\Omega^{x_2}_\eps \setminus \bar{C}^{x_2}_\eps} 
       W_\triangle(\nabla \tilde{y}) \, dx_1 \, dx_2
  \\& \geq \frac{2 \beta}{\gamma} \int_{\eps}^{1-\eps} k_{\eta}(x_2) \, dx_2 + \frac{2 c_R}{\gamma} \int_{D_1} k^C_{\eta}(x_2) \, dx_2 + \Big(\frac{\alpha l}{\sqrt{3}}\mu^2 + o(1)\Big) |D_2|  
    \\&\geq \frac{2 \beta}{\gamma} |I^{\eta}| + \frac{2 c_R}{\gamma}|D_1| + \Big(\frac{\alpha l}{\sqrt{3}}\mu^2 + o(1)\Big) |D_2| 
    \\& \geq \frac{2 \beta}{\gamma} |I^{\eta}|+ \min\left\{\frac{2 c_R}{\gamma}, \frac{\alpha l}{\sqrt{3}}\mu^2 + o(1)\right\}|I^\eta \setminus D^\mu|. 
\end{align*}
It follows $|I^{\eta} \setminus D^{\mu}| = O(\eps)$ and $|D^{\mu}| = 1- O(\eps)$, whence the crack energy from triangles in $\bar{{\cal C}}^{\mu}_{\eps,\eta}$ is given by $\frac{2\beta}{\gamma} + O(\eps)$ and then also $\eps \# \left(\bar{{\cal C}}_{\eps} \setminus \bar{{\cal C}}^{\mu}_{\eps,\eta}\right) = O(\eps)$. 

Finally, we concern ourselves with the projected distance of triangles in $\bar{{\cal C}}^{\mu}_{\eps,\eta}$. We first note that it suffices to show 
$$ \sup \left\{  |m_{\triangle_{1}} - m_{\triangle_{2}}| : \triangle_{1}, \triangle_{2} \in \bar{{\cal C}}^{\mu}_{\eps,\eta} \cap \bar{{\cal C}}^{(1)}_{\eps} \right\} = O(\eps)$$ 
since for a suitable $\tilde{\eta} \ge \eta$ for any $\triangle \in \bar{{\cal C}}^{\mu}_{\eps,\eta} \cap \bar{{\cal C}}^{(2)}_{\eps}$ there is a $\tilde{\triangle} \in \bar{{\cal C}}^{\mu}_{\eps,\tilde{\eta}} \cap \bar{{\cal C}}^{(1)}_{\eps}$ with $|m_{\triangle} - m_{\tilde{\triangle}}| \leq\eps$.  Let $x_2, z_{2} \in D^{\mu}$, $x_2 < z_2$ with $z_2 - x_2 \le C \eps$ and $|m_{\triangle_{1}} - m_{\triangle_{2}}| > 0$ for the corresponding broken triangles $\triangle_1, \triangle_2 \in \bar{\cal C}^{(1)}_{\eps}$. We may assume if a triangle intersects $(0, l) \times \{z_2\}$ or $(0, l) \times \{x_2\}$ then its interior does so, too. Denote by $\bar{d} = \gamma^{-1} |m_{\triangle_{1}} - m_{\triangle_{2}}|$ the distances of the centers in $\vv_{\gamma}$-projection onto the $x_1$-axis.

Let $x_1, z_1 \in (0,l)$ be the points on the slices $(0,l) \times \left\{x_2\right\}$ and $(0,l) \times \left\{z_2\right\}$ satisfying $\pi_{\vv^{\bot}_{\gamma}}(x_1,x_2) 
= m_{\triangle_1}$ and $\pi_{\vv^{\bot}_{\gamma}}(z_1,z_2) 
= m_{\triangle_2}$, respectively, where $\pi_{\vv^{\bot}_{\gamma}}$ denotes the orthogonal projection onto the linear subspace spanned by $\vv^{\bot}_{\gamma}$. Let $w= \e_1 \cdot \vv_{\gamma}|x_2 - z_2|/\gamma$. Then the $\vv_{\gamma}$-projection of $z = (z_1, z_2)$ onto the $x_2$-slice is given by $(\tilde{z}_1, x_2)$ with $\tilde{z}_1 = z_{1} - w$. Then $\bar{d} = |x_{1} - \tilde{z}_1|$ and without restriction we may assume $x_{1} > \tilde{z}_1$. 

Let $s_{\eps}= \frac{\sqrt{3}\eps}{4\gamma}$. We now consider the area bounded by the parallelogram with corners $(\tilde{z}_1 + s_{\eps},x_2)$, $(x_1 - s_{\eps}, x_2)$, $(z_1 + \bar{d} - s_{\eps}, z_2)$, $(z_1 + s_{\eps}, z_2)$. It is covered by $\frac{2\gamma \bar{d}}{\sqrt{3}\eps} - 1$ stripes of width $\frac{\sqrt{3}}{2}\eps$ in $\vv_{\gamma}$-direction consisting of lattice triangles intersecting the parallelogram, the first of these stripes touching $\triangle_1$, the last one touching $\triangle_2$ (note that if $\gamma \bar{d} = \frac{\sqrt{3}}{2}\eps$ the parallelogram is degenerated to a segment). For the intermediate stripes (\ref{eq: D_mu}) shows that 
\begin{align*} 
  y_1(t, x_2) 
  &\le t + \sqrt{\eps} l \mu \quad \forall\, t < x_1 - s_{\eps}
  \qquad \text{and} \\ 
  y_1(t, z_2) 
  &\ge t + \sqrt{\eps} l (a - \mu) \quad \forall\, t > z_1 + s_{\eps}. 
\end{align*}
This shows that if $(t, x_2)$ and $(t+w, z_2)$, $x_1 - \bar{d} + s_{\eps} < t < x_1 - s_{\eps}$ lie in the bottom and top triangles of some intermediate stripe, respectively, which are unbroken by construction of $D^{\mu}$, then 
\begin{align*} 
  |y(t+w, z_2) - y(t, x_2)| 
  \ge y_1(t+w, z_2) - y_1(t, x_2) 
  \ge w + \sqrt{\eps} l (a - 2 \mu) \sim \sqrt{\eps}.  
\end{align*}
Consider the $\frac{2\gamma \bar{d}}{\sqrt{3}\eps}$ atomic chains in $\vv_{\gamma}$ direction that lie on the boundary of these stripes. They are of length $\gamma^{-1} (z_2 - x_2) + O(\eps) \le C \eps \ll \sqrt{\eps}$.  So there is a constant $c > 0$ such that each of these chains contains at least one spring elongated by {a factor of more than $\frac{c}{\sqrt{\eps}}$. By passing, if necessary, to a lower threshold $\tilde{\eta} \ge \eta$, we obtain that the triangles sharing such a spring are broken and additionally one neighbor of each. As broken triangles for such springs on neighboring chains might overlap, we only consider every second atom chain and denote the set of type one triangles adjacent to such a spring on atom chains of odd numbers by $\bar{\cal C}^{(1)}_{\vv_{\gamma}}(\triangle_1, \triangle_2)$. We note that
\begin{equation}\label{eq: d bound}
\gamma \bar{d} \leq \sqrt{3} \eps \# \bar{\cal C}^{(1)}_{\vv_{\gamma}}(\triangle_1, \triangle_2).
\end{equation} 
The projection onto the $x_2$-axis of the spring in $\vv_{\gamma}$-direction is an interval $J$ of length $\gamma\eps$. Counting broken springs, it is elementary to see that the energy contribution $\frac{4}{\sqrt{3}\eps} \int_{(\eps,l-\eps)\times J} W_{\triangle}(\nabla \tilde{y}_{\eps})$ of the part of these broken triangles that lies in the stripe $(0, l) \times J$ is bounded from below by 
\begin{equation}\label{eq: bad strip}
2  \eps(  1 +  P(\gamma) ) \beta^{\tilde{\eta}},
\end{equation}
 where $P(\gamma)$ is the projection coefficient from \eqref{eq:P-gamma} satisfying $P(1) = \frac{1}{2}$ and in particular $P(\gamma) = 0$ if and only if $\gamma = \frac{\sqrt{3}}{2}$. On the other hand, the energy within stripes $(0, l) \times J'$ when $J'$ is the projection of an arbitrary broken triangle is still bounded from below by $2 \eps \beta^{\tilde{\eta}}$. 

Now let $\triangle_i$, $i = 1, \ldots, M_{\eps}$, denote all triangles $\triangle$ in $\bar{\cal C}^{\mu}_{\eps, \tilde{\eta}} \cap \bar{\cal C}^{(1)}_{\eps}$ such that there exists $x^{(i)}_2 \in D^{\mu}$ with $(0, l) \times \{x^{(i)}_2\}$ intersecting with the interior of $\triangle$. The numbering shall be chosen so as to satisfy $x_2^{(1)} < \ldots < x_2^{(M_{\eps})}$. As $1 - |D^{\mu}| = O(\eps)$, there exists a constant $C > 0$ such that $x_2^{(i+1)} - x_2^{(i)} < C \eps$, $i = 1, \ldots, M_{\eps}-1$. We define the subset $\{x^{(i_{j})}_2\}_{j=1, \ldots N_{\eps}}$ of $\{x^{(i)}_2\}_{i=1, \ldots, M_{\eps}}$ such that $x^{({i})}_2 = x^{(i_{j})}_2$ for a $j = 1, \ldots N_{\eps}$ if and only if $|m_{\triangle_{i}} - m_{\triangle_{i+1}}| > 0$. According to our previous considerations, if $I^{\tilde{\eta}}_{\vv_{\gamma}}$ is the projection of $\bar{\cal C}^{(1)}_{\vv_{\gamma}} := \bigcup_{j = 1}^{N_{\eps}} \bar{\cal C}^{(1)}_{\vv_{\gamma}}(\triangle_{i_{j}}, \triangle_{i_{j}+1})$ onto the $x_2$-axis, then 
\begin{equation}\label{eq: I_gamma}
|I^{\tilde{\eta}}_{\vv_{\gamma}}| \leq \gamma \eps  \# \bar{\cal C}^{(1)}_{\vv_{\gamma}}.
\end{equation}
As before using (\ref{eq: bad strip}) and (\ref{eq: I_gamma}) we see that the total energy is greater or equal to 
\begin{align*} 
  & \# \bar{\cal C}^{(1)}_{\vv_{\gamma}} 2 \eps ( 1 + P(\gamma) ) \beta^{\tilde{\eta}} 
   +  |I^{\tilde{\eta}} \setminus I^{\tilde{\eta}}_{\vv_{\gamma}}| \frac{2 \beta^{\tilde{\eta}}}{\gamma} + O(\eps) \\ 
    &= |I^{\tilde{\eta}}| \frac{2 \beta^{\tilde{\eta}}}{\gamma} + 2 \# \bar{\cal C}^{(1)}_{\vv_{\gamma}} \eps P(\gamma) \beta^{\tilde{\eta}} + 2 \# \bar{\cal C}^{(1)}_{\vv_{\gamma}} \eps  \beta^{\tilde{\eta}} - |I^{\tilde{\eta}}_{\vv_{\gamma}}| \frac{2 \beta^{\tilde{\eta}}}{\gamma}+  O(\eps) \\ 
    & \geq \frac{2 \beta}{\gamma} + 2 \# \bar{\cal C}^{(1)}_{\vv_{\gamma}} \eps P(\gamma) \beta^{\tilde{\eta}} + O(\eps), 
\end{align*} 
and so $\# \bar{\cal C}^{(1)}_{\vv_{\gamma}} = O(1)$. As every $\triangle \in \bar{\cal C}^{(1)}_{\vv_{\gamma}}$ is in at most two different $\bar{\cal C}^{(1)}_{\vv_{\gamma}}(\triangle_{i_{j}}, \triangle_{i_{j}+1})$, this also yields $\sum_{j = 1}^{N_{\eps}} \# \bar{\cal C}^{(1)}_{\vv_{\gamma}}(\triangle_{i_{j}}, \triangle_{i_{j}+1}) = O(1)$.

Applying (\ref{eq: d bound}) we find that 
\begin{align*} 
  O(1) 
  = \sum_{j = 1}^{N_{\eps}} \# \bar{\cal C}^{(1)}_{\vv_{\gamma}}(\triangle_{i_{j}}, \triangle_{i_{j}+1}) 
  \ge \sum_{j = 1}^{N_{\eps}} \frac{\gamma \bar{d}_{i_{j}}}{\sqrt{3} \eps}  
  \ge \frac{c}{\eps} \sum_{j = 1}^{N_{\eps}} |m_{\triangle_{i_{j}}} - m_{\triangle_{i_{j}+1}}|
\end{align*}
for a constant $c > 0$, when $\bar{d}_i = \gamma^{-1} |m_{\triangle_{i}} - m_{\triangle_{i + 1}}|$. This concludes the proof. 
\eop

The above Lemmas \ref{lemma: minimizing sequences - properties - elast} and \ref{lemma: minimizing sequences - properties - crack} show that for a sequence of almost minimizers $(\tilde{y}_{\eps})$ satisfying \eqref{eq:almost-min}, the number $\#\bar{\cal C}_{\eps}$ of largely deformed triangles is bounded independently of $\eps$ for $a < a_{\rm crit}$, while in the supercritical case for $\phi \ne 0$ there are two subsets 
\begin{align}\label{eq:p-eps}
\begin{split}
  \Omega^{(1)}_{\eps} 
  &:= \left\{x \in \Omega_{\eps} : 0 \leq x_1 \leq p_{\eps} - c\eps + (\vv_{\gamma} \cdot \e_1) x_2  \right\},\\
  \Omega^{(2)}_{\eps} 
  &:= \left\{x \in \Omega_{\eps} : p_{\eps} + c\eps + (\vv_{\gamma} \cdot \e_1) x_2\leq x_1 \leq  l  \right\}, 
\end{split}
\end{align}
$c > 0$ independent of $\eps$ and $p_{\eps}$ to be chosen appropriately, such that the number of triangles in $\bar{\cal C}_{\eps}$ intersecting $\Omega^{(1)}_{\eps} \cup \Omega^{(2)}_{\eps}$ is bounded uniformly in $\eps$. We recall that the last claim in Lemma \ref{lemma: minimizing sequences - properties - crack} does not hold if $\vv_{\gamma}$ is not unique ($\gamma = \frac{\sqrt{3}}{2}$). Indeed, if $P(\gamma)$ vanishes, we cannot conlude that $\# \bar{\cal C}^{(1)}_{\vv_{\gamma}} = O(1)$ in the above proof. In this case we do not expect that the essential part of the broken triangles lies in in a small stripe parallel to $\R (\frac{1}{2}, \frac{\sqrt{3}}{2})^T$ or $\R (\frac{1}{2}, \frac{\sqrt{3}}{2})^T$ as we have already seen that the crack can take a serrated course. Nevertheless, if $\gamma = \frac{\sqrt{3}}{2}$ (or equivalently $\phi = 0$) one can show that up to a number being uniformly bounded in $\eps$ the broken triangles $\bar{\cal C}_{\eps}$ lie in a stripe around the graph of a Lipschitz function. 

\begin{lemma}\label{lemma: gamma special case}
Let $\tilde{u}_{\eps}$ be a minimizing sequence satisfying 
$${\cal E}_{\eps}(\id + \sqrt{\eps} u_{\eps}) = \inf {\cal E}_{\eps}+O(\eps).$$
Let $a> a_{\rm crit}$ and $\phi = 0$. Then there exists a Lipschitz function $g:(0,1) \to (\psi(\eps),l - \psi(\eps))$ with $|g'| = \frac{1}{\sqrt{3}}$ a.e.\ such that for $\mu$ sufficiently small $\eps \# (\bar{\cal C}_\eps \setminus \bar{\cal C}^{\mu}_{\eps,\eta}) = O(\eps)$ and 
\begin{align}\label{eq: Lipschitz stripe}
  \cup_{\triangle \in \bar{\cal C}^{\mu}_{\eps,\eta}} \triangle 
  \subset \left\{(x,y) \in \Omega: g(y) - C\eps \leq x \leq g(y) + C\eps\right\},
\end{align}
for some $C > 0$ independent of $g$ and $\eps$.
\end{lemma}

\Proof We have $\vv_2 = (\frac{1}{2}, \frac{\sqrt{3}}{2})^T$, $\vv_3 := \vv_2 - \vv_1 = (-\frac{1}{2}, \frac{\sqrt{3}}{2})^T$ and $\vv^{\bot}_2 = ( - \frac{\sqrt{3}}{2}, \frac{1}{2})^T$, $\vv^{\bot}_3 = -( \frac{\sqrt{3}}{2}, \frac{1}{2})^T$. By Lemma \ref{lemma: minimizing sequences - properties - crack} we immediately get $|I^{\eta}| =1 - O(\eps)$ and $\eps \# (\bar{\cal C}_\eps \setminus \bar{\cal C}^{\mu}_{\eps,\eta}) = O(\eps)$ for $\mu$ sufficiently small recalling that these properties were derived independently of the choice of $\gamma$. Similarly as before we note that after passing to a suitable $\tilde{\eta} \ge \eta$ it suffices to show the claim for $\tilde{\triangle} \in \bar{{\cal C}}^{\mu}_{\eps,\tilde{\eta}} \cap \bar{{\cal C}}^{(1)}_{\eps}$. We estimate the difference of broken triangles $\bar{\cal C}^{\mu}_{\eps,\eta} \cap \bar{\cal C}^{(1)}_{\eps}$ projected onto the linear subspaces spanned by $\vv^\bot_2$ and $\vv^\bot_3$. We recall that the projection of some triangle $\triangle$ on these subspaces are intervals of length $\frac{\sqrt{3}}{2}\eps$ and denote the centers of the intervals by $m^{(2)}_{\triangle}$ and $m^{(3)}_{\triangle}$, respectively. 

Let $x_2, z_{2} \in D^{\mu}$, $x_2 < z_2$ with $z_2 - x_2 \le C \eps$ and 
\begin{align}\label{eq: projected distance1}
  (m^{(2)}_{\triangle_{1}} - m^{(2)}_{\triangle_{2}}) \cdot \vv^{\bot}_2 > 0 
\end{align}
or 
\begin{align}\label{eq: projected distance2}
  (m^{(3)}_{\triangle_{1}} - m^{(3)}_{\triangle_{2}}) \cdot \vv^{\bot}_3 < 0
\end{align}
for the corresponding broken triangles $\triangle_1, \triangle_2 \in \bar{\cal C}^{(1)}_{\eps}$. Without restriction we treat the case (\ref{eq: projected distance1}). As in the proof of Lemma \ref{lemma: minimizing sequences - properties - crack} we may assume if a triangle intersects $(0, l) \times \{z_2\}$ or $(0, l) \times \{x_2\}$ then its interior does so, too. Denote by $\bar{d}^{(i)} = \frac{2}{\sqrt{3}} |m^{(i)}_{\triangle_{1}} - m^{(i)}_{\triangle_{2}}|$, $i=2,3$, the distances of the centers in $\vv_{i}$-projection onto the $x_1$-axis.

Let $x^{(i)}_1, z^{(i)}_1 \in (0,l)$ such that $\pi_{\vv^{\bot}_{i}}(x^{(i)}_1,x_2) 
= m^{(i)}_{\triangle_1}$ and $\pi_{\vv^{\bot}_{i}}(z^{(i)}_1,z_2) 
= m^{(i)}_{\triangle_2}$, respectively, where $\pi_{\vv^{\bot}_{i}}$ denotes the orthogonal projection onto the linear subspace spanned by $\vv^{\bot}_{i}$, $i=2,3$. Let $w^{(2)} = \frac{1}{\sqrt{3}}|x_2 - z_2|$ and $w^{(3)} = -\frac{1}{\sqrt{3}}|x_2 - z_2|$. Then the $\vv_{i}$-projection of $z^{(i)} = (z^{(i)}_1, z_2)$ onto the $x_2$-slice is given by $(\tilde{z}^{(i)}_1, x_2)$ with $\tilde{z}^{(i)}_1 = z^{(i)}_{1} - w^{(i)}$ for $i=2,3$. We note that $\bar{d}^{(i)} = |x^{(i)}_{1} - \tilde{z}^{(i)}_1|$. Taking (\ref{eq: projected distance1}) into account we obtain  $x^{(2)}_{1} < \tilde{z}^{(2)}_1 < \tilde{z}^{(3)}_1$.

Let $s_{\eps}= \frac{\eps}{2}$. As in the previous proof we consider areas bounded by parallelograms. For $i=2,3$, let $P^{(i)}$ be the parallelogram with corners $(x^{(i)}_1 + s_{\eps}, x_2)$, $(\tilde{z}^{(i)}_1 - s_{\eps},x_2)$, $(z^{(i)}_1 - s_{\eps}, z_2)$,  $(z^{(i)}_1 - \bar{d}^{(i)} + s_{\eps}, z_2)$. They are covered by $\frac{\bar{d}^{(i)}}{\eps} - 1$ stripes of width $\frac{\sqrt{3}}{2}\eps$ in $\vv_{i}$-direction, respectively (note that $P^{(2)}$ can again be degenerated to a segment if $\bar{d}^{(2)} = \eps$). It is not hard to see that both parallelograms cover $\left\lceil \frac{2 |z_2 - x_2|}{\sqrt{3} \eps}\right\rceil$ or $\left\lceil \frac{2 |z_2 - x_2|}{\sqrt{3} \eps}\right\rceil + 1$ stripes of width $\frac{\sqrt{3}}{2}\eps$ in $\e_1$-direction, where the stripes at the top and at the bottom are only partially covered (the exact number depends of the precise location of the slices $(0,l) \times \{x_2\}$ and $(0,l) \times \{z_2\}$). We denote the number of these covered stripes by $N(\triangle_1,\triangle_2)$ and the orthogonal projection onto the $x_2$-axis by $I(\triangle_1,\triangle_2)$. Setting 
\begin{align}\label{eq: n definition}
  n^{(i)}(\triangle_1,\triangle_2) = \frac{2}{\sqrt{3}\eps}|m^{(i)}_{\triangle_{1}} - m^{(i)}_{\triangle_{2}}|
\end{align}
it is elementary to see that $\bar{d}^{(2)} = n^{(2)}(\triangle_1,\triangle_2) \eps$ and $\bar{d}^{(3)} = (n^{(2)}(\triangle_1,\triangle_2) + N(\triangle_1,\triangle_2) - 1)\eps  $.

Following the lines of the previous proof we see that each of the $\frac{\bar{d}^{i}}{\eps}$ atomic chains in $\vv_{i}$ direction lying on the boundary of the stripes which cover $P^{(i)}$, contains at least one spring elongated by a factor of more than $\frac{c}{\sqrt{\eps}}$. Consequently, on the $N(\triangle_1,\triangle_2)$ stripes in $\e_1$-direction we have at least $2 n^{(2)}(\triangle_1,\triangle_2) + N(\triangle_1,\triangle_2) - 1 > N(\triangle_1,\triangle_2)$ broken springs orientated in $\vv_2$ or $\vv_3$ direction.
Let $J$ be an interval of length $\frac{\sqrt{3}}{2}\eps$ such that the stripe $(0,l) \times J$ consists of lattice triangles. 
It is elementary to see that if two broken springs in $\vv_2$ and $\vv_3$ lie in the stripe at least three triangles are broken, i.e. lie in the set $\bar{\cal C}_{\eps,\tilde{\eta}}$. Thus, the energy contribution $\frac{4}{\sqrt{3}\eps} \int_{(\eps,l-\eps)\times J} W_{\triangle}(\nabla \tilde{y}_{\eps})$ of the stripe can be bounded from below by $3\eps \beta^{\tilde{\eta}}.$
More generally, if on a stripe there are $k \in \N$ broken springs in $\vv_2$ and $\vv_3$ the energy contribution is at least
$(k+1)\eps \beta^{\tilde{\eta}}.$

On the other hand, we recall that on an arbitrary stripe $(0,l) \times J'$ consisting of lattice triangles the energy is always bounded from below by $2\eps\beta^{\tilde{\eta}}$. Consequently, we derive that in the above situation the energy contribution of the $N(\triangle_1,\triangle_2)$ stripes $\frac{4}{\sqrt{3}\eps} \int_{(\eps,l-\eps)\times I(\triangle_1,\triangle_2)} W_{\triangle}(\nabla \tilde{y}_{\eps})$ is bounded from below by 
\begin{align}\label{eq: stripe energy}
N(\triangle_1,\triangle_2)\eps \beta^{\tilde{\eta}} + (2n^{(2)}(\triangle_1,\triangle_2) + N(\triangle_1,\triangle_2) - 1) \eps \beta^{\tilde{\eta}}
\end{align}
and note that the energy contribution of $N(\triangle_1,\triangle_2)$ stripes is always bounded from below by $2N(\triangle_1,\triangle_2) \eps\beta^{\tilde{\eta}}$. 
 
Now let $\triangle_i$, $i = 1, \ldots, M_{\eps}$, denote all triangles $\triangle$ in $\bar{\cal C}^{\mu}_{\eps, \tilde{\eta}} \cap \bar{\cal C}^{(1)}_{\eps}$ such that there exists $x^{(i)}_2 \in D^{\mu}$ with $(0, l) \times \{x^{(i)}_2\}$ intersecting with the interior of $\triangle$. The numbering shall be chosen so as to satisfy $x_2^{(1)} < \ldots < x_2^{(M_{\eps})}$. As $1 - |D^{\mu}| = O(\eps)$, there exists a constant $C > 0$ such that $x_2^{(i+1)} - x_2^{(i)} < C \eps$, $i = 1, \ldots, M_{\eps}-1$. We define the subset $\{x^{(i_{j})}_2\}_{j=1, \ldots N_{\eps}}$ of $\{x^{(i)}_2\}_{i=1, \ldots, M_{\eps}}$ such that $x^{({i})}_2 = x^{(i_{j})}_2$ for a $j = 1, \ldots N_{\eps}$ if and only if $m_{\triangle_{i}}$ and $m_{\triangle_{i+1}}$ satisfy (\ref{eq: projected distance1}) or (\ref{eq: projected distance2}). We let $p(j) = 2$ or $p(j) = 3$ if (\ref{eq: projected distance1}) or (\ref{eq: projected distance2}) holds, respectively. Let $I_\eps = \cup^{N_\eps}_{j =1} I(\triangle_{i_j}, \triangle_{i_j + 1})$. Taking ({\ref{eq: stripe energy}}) into account the energy contribution $\frac{4}{\sqrt{3}\eps} \int_{(\eps,l-\eps)\times I_\eps} W_{\triangle}(\nabla \tilde{y}_{\eps})$ can be bounded from below by 
$$ \frac{4}{\sqrt{3}\eps}|I_{\eps}| \eps\beta^{\tilde{\eta}} + \frac{1}{2} \sum^{N_\eps}_{i=1} (2n^{(p(j))}(\triangle_{i_j},\triangle_{i_j + 1}) - 1) \eps \beta^{\tilde{\eta}} + O(\eps). $$
The factor $\frac{1}{2}$ accounts for the possibility that two adjacent intervals $I(\triangle_{i_j},\triangle_{i_j + 1})$, $I(\triangle_{i_{j+1}},\triangle_{i_{j+1} + 1})$ may overlap. We thus see that the total energy is greater or equal to
$$\frac{4}{\sqrt{3}}\beta^{\tilde{\eta}} + \frac{1}{2} \sum^{N_\eps}_{i=1} (2n^{(p(j))}(\triangle_{i_j},\triangle_{i_j + 1}) - 1) \eps \beta^{\tilde{\eta}} + O(\eps) $$
and so $\sum^{N_\eps}_{j=1} n^{(p(j))}(\triangle_{i_j},\triangle_{i_j + 1}) = O(1)$. We now construct the function $g: (0,1) \to (0,l)$. For $i \in M_\eps$ let $M^{(1)}_i$ and $M^{(2)}_i$ be the orthogonal projections of the center of $\triangle_i$ onto the $x_1$ and $x_2$-axis, respecively. If $i \notin N_\eps$ set 
$$\tilde{\tilde{g}}= \frac{M^{(1)}_{i+1} - M^{(1)}_i}{M^{(2)}_{i+1} - M^{(2)}_i}$$
on the interval $[M^{(2)}_{i}, M^{(2)}_{i+1}]$. Now let $\tilde{g}$ be the Lipschitz function satisfying $\tilde{g}(M^{(2)}_1) = M^{(1)}_1$ and $\tilde{g}' = \tilde{\tilde{g}}$. By construction it is easy to see that $|\tilde{g}'|\leq \frac{1}{\sqrt{3}}$ on $[M^{(2)}_1, M^{(2)}_{N_\eps}]$. We extend $\tilde{g}$ arbitrarily to $(0,1)$ such that $\left\|\tilde{g}'\right\|_\infty \leq \frac{1}{\sqrt{3}}$. By (\ref{eq: n definition}) we have
$$\sum^{N_\eps}_{j=1} |m^{(p(j))}_{\triangle_{i_j}}-m^{(p(j))}_{\triangle_{i_j + 1}}| = O(\eps), $$
and then is not hard to see that there is some $C>0$ independent of $\tilde{g}$ and $\eps$ such that (\ref{eq: Lipschitz stripe}) holds. Recalling \eqref{eq: bv-phi-0} it remains to choose $g : (0,1) \to (\psi(\eps), l - \psi(\eps))$ with $g' = \pm \frac{1}{\sqrt{3}}$ a.e.\ and $\| g - \tilde{g} \|_{\infty} \le C \eps$. \eop

We conclude that for $\phi = 0$ in the supercritical case there are two subsets 
\begin{align}\label{eq: g}
\begin{split}
  \Omega^{(1)}_{g_\eps} 
  &:= \left\{x \in \Omega_{\eps} : 0 \leq x_1 \leq g_\eps(x_2) - c\eps  \right\},\\
  \Omega^{(2)}_{g_\eps} 
  &:= \left\{x \in \Omega_{\eps} : c\eps + g_\eps(x_2) \leq  x_1 \leq  l  \right\}, 
\end{split}
\end{align}
where $g_\eps$ is chosen appropriately as in Lemma \ref{lemma: gamma special case} and $c > 0$ independent of $\eps$, such that the number of triangles in $\bar{\cal C}_{\eps}$ intersecting $\Omega^{(1)}_{g_\eps} \cup \Omega^{(2)}_{g_\eps}$ is bounded uniformly in $\eps$.
Note that with $\eps \ll \bar{\psi}(\eps) = \psi(\eps) - c \eps \ll 1$ one has 
\begin{align}\label{eq: components connected}
\begin{split}
  \left( (0,\bar{\psi}(\eps)) \times (0,1) \right) \cap \Omega_\eps \subset \Omega^{(1)}_{g_\eps}, \\ 
  \left( (l - \bar{\psi}(\eps),l) \times (0,1)\right) \cap \Omega_\eps  \subset \Omega^{(2)}_{g_\eps}, 
\end{split}
\end{align}
so that, in particular, $\Omega^{(1)}_{g_\eps}$ and $\Omega^{(2)}_{g_\eps}$ are connected.
The following lemma shows that broken triangles in these sets can be `healed'. In order to treat these cases simultaneously in the following we will call these sets the `good set'
$$ \Omega_{\rm good} 
   = \begin{cases} 
       \Omega_{\eps} 
       &\text{ for } a < a_{\rm crit}, \\ 
       \Omega^{(1)}_{\eps} \cup \Omega^{(2)}_{\eps} 
       &\text{ for } a > a_{\rm crit},\, \phi \ne 0 \text{ and } \\ 
       \Omega^{(1)}_{g_\eps} \cup \Omega^{(2)}_{g_\eps} 
       &\text{ for } a > a_{\rm crit},\, \phi = 0,
     \end{cases} $$
with $\Omega^{(i)}_{\eps}$ and $\Omega^{(i)}_{g_\eps}$, $i = 1,2$, as defined in \eqref{eq:p-eps} and \eqref{eq: g}.

\begin{lemma}\label{lemma:triangle-healing}
Suppose $\tilde{y}_{\eps}$ is a minimizing sequence satisfying 
${\cal E}_{\eps}(y_{\eps}) = \inf {\cal E}_{\eps}+O(\eps)$. 
There exists $\bar{y}_{\eps} \in W^{1, \infty}(\Omega_{\rm good}; \R^2)$ with $\nabla \bar{y}_{\eps}$ bounded in $L^{\infty}(\Omega_{\rm good})$ uniformly in $\eps$ such that 
\begin{align*} 
  |\{x \in \Omega_{\rm good} : \bar{y}_{\eps}(x) \ne \tilde{y}_{\eps}(x)\}| 
    = O(\eps^2) 
\end{align*}
and 
\begin{align*} 
  \int_{\Omega_{\rm good}} \dist^2(\nabla \bar{y}_{\eps}(x), SO(2)) \, dx 
  \le C \int_{\Omega_{\rm good} \setminus \bigcup_{\triangle \in \bar{\cal C}_{\eps}} \triangle} 
      \dist^2(\nabla \tilde{y}_{\eps}, SO(2)) \, dx. 
\end{align*}
\end{lemma}

\Proof For notational convenience we drop the subscript $\eps$ in the following proof. By Lemmas \ref{lemma: minimizing sequences - properties - elast}, \ref{lemma: minimizing sequences - properties - crack} and \ref{lemma: gamma special case} we can partition the area covered by the (closed) triangles in $\bar{\cal C}$ intersecting $\Omega_{\rm good}$ into connected components $C_1, \ldots, C_N$ such that 
$$ \bigcup_{\triangle \in \bar{\cal C} : \triangle \cap \Omega_{\rm good} \ne \emptyset} \triangle  
   = C_1 \dot{\cup} \ldots \dot{\cup} C_N, $$ 
where $N$ is bounded uniformly in $\eps$. Then the maximal diameter of each set $C_i$ is bounded by a term $O(\eps)$. For each $i$, the largest connected component $D_i$ of the complement $\Omega_{\rm good} \setminus C_i$ lying in the same component of $\Omega_{\rm good}$ is unique (with area of the order 1 while all the other components of the complement are of size $O(\eps^2)$). Let $V_i$ be the union of triangles whose interior is contained in $D_i$ that touch the boundary of $C_i$. 

We now proceed to define $\bar{y}$ by modifying $\tilde{y}$ on all the triangles not contained in $\overline{D}_i$, successively for $i = 1, \ldots, N$. For each $i$ this modification is done iteratively on triangles $\triangle$ which share at least one side with a triangle that has been modified previously or with a triangle lying in $V_i$ in such a way that $\bar{y}$ is continuous along such sides and $\bar{y}|_{\triangle}$ is affine and minimizes $\dist((\bar{y})_{\triangle}, SO(2))$. 

In order to estimate $\dist(\nabla \bar{y}, SO(2))$ we recall the following geometric rigidity result proved in \cite{FrieseckeJamesMueller:02}: If $U \subset \R^d$ is a (connected) Lipschitz domain, then there exists a constant $C = C(U)$ such that for any $f \in H^1(U, \R^d)$ there is a rotation $R \in SO(d)$ with 
\begin{align}\label{eq:geometric-rigidity} 
  \int_U |\nabla f(x) - R|^2 \, dx 
  \le C \int_U \dist^2(\nabla f(x), SO(d)) \, dx. 
\end{align}
The constant $C(U)$ is invariant under rescaling of the domain. For later use we mention that if $\dist^2(\nabla f(x), SO(d))$ is equiintegrable, then $R$ can be chosen in such a way that also $|\nabla f(x) - R|^2$ is equiintegrable, cf.\ \cite{FrieseckeJamesMueller:06}. 

Consider a single step in the modification process, when $\tilde{y}$ is modified to $\bar{y}$ on $\triangle$, and let $U$ be the union of triangles that have been modified previously or lie in $V_i$. By the geometric rigidity estimate \eqref{eq:geometric-rigidity}, there is a rotation $R \in SO(2)$ such that \eqref{eq:geometric-rigidity} holds for $f = \bar{y}$. Since $\nabla \bar{y}$ is piecewise constant, this means 
$$ \sum_{\triangle' \subset U} |(\bar{y})_{\triangle'} - R|^2 
   \le C \sum_{\triangle' \subset U} \dist^2((\bar{y})_{\triangle'}, SO(2)). $$ 
It is not hard to see that there exists an extension $w$ of $\bar{y}$ from $U$ to $U \cup \triangle$ such that 
$$ |(w)_{\triangle} - R|^2 
   \le C \sum_{\triangle' \subset U} |(\bar{y})_{\triangle'} - R|^2. $$ 
(If there is only one side of $\triangle$ on the boundary of $U$, say adjacent to $\triangle' \subset U$, then one can take $w$ with $(w)_{\triangle} = (\bar{y})_{\triangle'}$. If at least two sides, say in $\vv_1$ and $\vv_2$ direction, are shared by triangles $\triangle_1, \triangle_2 \subset U$, respectively, then these sides have a common corner and the unique extension $w$ satisfies $(w)_{\triangle}\vv_i = (\bar{y})_{\triangle_i}\vv_i = R \vv_i + ((\bar{y})_{\triangle_i} - R) \vv_i$, $i = 1, 2$.) Now by construction of $\bar{y}$ on $\triangle$ we see that 
$$ \dist^2((\bar{y})_{\triangle}, SO(2)) 
   \le C \sum_{\triangle' \subset U} |(\bar{y})_{\triangle'} - R|^2 $$ 
and so 
\begin{align*} 
  \int_{U \cup \triangle} \dist^2(\nabla \bar{y}(x), SO(2)) \, dx 
  \le C \int_{U} \dist^2(\nabla \bar{y}(x), SO(2)) \, dx. 
\end{align*}
Iterating this estimate we finally arrive at 
\begin{align*} 
  \int_{\Omega_{\rm good}} \dist^2(\nabla \bar{y}(x), SO(2)) \, dx 
  \le C \int_{\Omega_{\rm good} \setminus \bigcup_{i} C_i} 
      \dist^2(\nabla \tilde{y}, SO(2)) \, dx. 
\end{align*}
Here the constant $C$ can be chosen independently of $\eps$. This is due to the facts that the number of modification steps is bounded uniformly in $\eps$ and -- after rescaling the shapes $U$ with $\frac{1}{\eps}$ -- there is also only a uniformly bounded number of shapes $U$ involved in the previous rigidity estimates. Moreover, each triangle is covered by no more than three of the sets $V_i$. 

The uniform boundedness of the number of modification steps also shows that $|\{x \in \Omega_{\rm good} : \bar{y}(x) \ne \tilde{y}(x)\}| = O(\eps^2)$ and, by definition of $\bar{\cal C}$ and construction of $\bar{y}$, that $\|\nabla \bar{y}\|_{L^{\infty}(\Omega_{\rm good})} = O(1)$. \eop 

Note that up to a set of small size $\bar{y}_{\eps}$ satisfies the same boundary conditions as $\tilde{y}_{\eps}$ on the lateral boundary. More precisely, there are $\Gamma^{(i)}_{\eps} \subset (0,1)$, $|\Gamma^{(i)}_{\eps}| = O(\eps)$, $i=1,2$, such that $\bar{y}_{\eps}$ and $\tilde{y}_{\eps}$ coincide on $(0,\eps) \times ((0,1)\setminus \Gamma^{(1)}_{\eps})$ and $(l-\eps,l) \times ((0,1)\setminus \Gamma^{(2)}_{\eps})$. With these boundary conditions and the geometric rigidity estimate \eqref{eq:geometric-rigidity} we can now derive strong convergence results for $\bar{y}_{\eps}$ and even the corresponding rescaled displacement $\bar{u}_{\eps} = \frac{1}{\sqrt{\eps}}(\bar{y}_{\eps} - \id)$ on $\Omega_{\rm good}$. We first consider the supercritical case and treat the cases $\phi \ne 0$ and $\phi = 0$ separately.
\begin{lemma}\label{lemma:strong-convergence-crack} 
If $a > a_{\rm crit}$ and $\phi \ne 0$, then there exist sequences $s_{\eps}, t_{\eps} \in \R$ such that 
$$ \| \bar{u}_{\eps} - (0, s_{\eps}) \|_{H^1(\Omega^{(1)}_{\eps})} 
   + \| \bar{u}_{\eps} - (al, t_{\eps}) \|_{H^1(\Omega^{(2)}_{\eps})} 
   \to 0. $$ 
\end{lemma} 

\Proof 
We again drop the subscript $\eps$. By applying the geometric rigidity estimate \eqref{eq:geometric-rigidity} to $\Omega^{(1)}$ and to $\Omega^{(2)}$, we obtain rotations $R^{(1)}, R^{(2)} \in SO(2)$ such that 
\begin{align}\label{eq: rotations} \| \nabla \bar{y} - R^{(i)} \|_{L^2(\Omega^{(i)}_{\eps})}
   \le C \| \dist(\nabla \bar{y}, SO(2)) \|_{L^2(\Omega^{(i)}_{\eps})}, 
   \quad i = 1, 2.
\end{align}
Here $C$ can be chosen independently of $\eps$ as all the possible shapes of $\Omega^{(i)}$ are related through bi-Lipschitzian homeomorphisms with Lipschitz constants of both the homeomorphism itself and its inverse bounded uniformly in $\eps$, cf.\ \cite{FrieseckeJamesMueller:02}. Now using that $\nabla \bar{y}$ is uniformly bounded in $L^{\infty}$, we obtain from Lemmas \ref{lemma:triangle-healing} and \ref{lemma:quadratic lower bound}(i)
\begin{align*}
  \sum_{i = 1}^2 \| \nabla \bar{y} - R^{(i)} \|_{L^2(\Omega^{(i)}_{\eps})}^2 
  &\le C \int_{\Omega_{\rm good} \setminus \bigcup_{\triangle \in \bar{\cal C}_{\eps}} \triangle} 
      \dist^2(\nabla \tilde{y}, SO(2)) \, dx \\ 
      &\le C \int_{\Omega_{\rm good} \setminus \bigcup_{\triangle \in \bar{\cal C}_{\eps}} \triangle} 
      \dist^2(\nabla \tilde{y}, O(2)) + \chi(\nabla \tilde{y})\, dx \\ 
  &\le C \int_{\Omega_{\rm good} \setminus \bigcup_{\triangle \in \bar{\cal C}_{\eps}} \triangle} 
      W_{\triangle,\chi}(\nabla \tilde{y}) \, dx. 
\end{align*}
But, as seen before, 
\begin{align*}
  \frac{4}{\sqrt{3} \eps} \int_{\Omega_{\rm good} \setminus \bigcup_{\triangle \in \bar{\cal C}_{\eps}} \triangle} 
      W_{\triangle, \chi}(\nabla \tilde{y}) \, dx  
  &\le {\cal E}^{\chi}(y) - \frac{2 \beta^{\eta}}{\gamma} |I^{\eta}| 
  = O(\eps) 
\end{align*}
by Lemma \ref{lemma: minimizing sequences - properties - crack}, and so 
\begin{align*}
  \sum_{i = 1}^2 \| \nabla \bar{y} - R^{(i)} \|_{L^2(\Omega^{(i)}_{\eps})}^2 
  = O(\eps^2).  
\end{align*}
By Poincar{\'e}'s inequality we then deduce that there are $\zeta^{(i)} \in \R^2$ such that 
\begin{align}\label{eq:y-on-rigid-components}
  \sum_{i = 1}^2 \| \bar{y} - R^{(i)} \cdot - \zeta^{(i)} \|_{H^1(\Omega^{(i)}_{\eps})} 
  = O(\eps). 
\end{align}
We extend $\bar{y}$ as an $H^1$-function from $\Omega^{(i)}_{\eps}$ to $\Omega^{(i)}$ (as defined in Theorem \ref{theo: conv-of-minimizers}), $i=1,2$, such that (\ref{eq:y-on-rigid-components}) still holds and $\bar{y}_1(0,x_2) = 0$ for $x_2 \in (0,1)\setminus \Gamma^{(1)}_{\eps}$, $\bar{y}_1(l,x_2) = l(1 + a_{\eps})$ for $x_2 \in (0,1)\setminus \Gamma^{(2)}_{\eps}$. The trace theorem for Sobolev functions with $x_1 = 0$ or $x_1 = l$ according to $i = 1$ and $i = 2$, respectively, gives  
\begin{align*}
  \sum_{i = 1}^2 \| \bar{y}(x_1, \cdot) - R^{(i)}(x_1, \cdot) - \zeta^{(i)} \|_{L^2(0,1)} 
  = O(\eps).  
\end{align*}
In particular, setting $\tilde{\zeta}^{(1)} = \zeta^{(1)}$ and $\tilde{\zeta}^{(2)} = \zeta^{(2)} - l a_{\eps}\e_1$, the first components satisfy 
\begin{align}\label{eq:first-component-bdry}
  \sum_{i = 1}^2 \| x_1  - R^{(i)}_{11} x_1 -  R^{(i)}_{12} \cdot - \tilde{\zeta}^{(i)}_1 \|_{L^2((0,1)\setminus \Gamma^{(i)}_{\eps})} 
  = O(\eps).  
\end{align} 
But then also the constant function 
\begin{align*}
  \frac{1}{2} R^{(i)}_{12}
  = \left( x_1  - R^{(i)}_{11} x_1 - R^{(i)}_{12} \left(\cdot - \frac{1}{2}\right) - \tilde{\zeta}^{(i)}_1 \right) 
    - \left( x_1  - R^{(i)}_{11} x_1 - R^{(i)}_{12} \cdot - \tilde{\zeta}^{(i)}_1 \right) 
\end{align*}
is of order $\eps$ in $L^2((\frac{1}{2}, 1) \setminus \Gamma^{(i)})$ and thus $|R^{(i)}_{12}| \le C\eps$. An elementary argument now yields 
$$ |R^{(i)} - \Id| = O(\eps) 
   \qquad \text{or} \qquad 
   |R^{(i)} + \Id| = O(\eps).$$ 
It is not hard to see that $|R^{(i)} - \Id| = O(\eps)$ as otherwise, e.g. for $i=1$, on the set $T =\left\{\triangle \in {\cal C}_\eps \setminus \bar{\cal C}_\eps: \triangle \subset (0,\eps) \times (0,1)\right\}$ we get, due to the boundary conditions, 
$$O(\eps^2) = \int_{T} |\nabla \bar y - R^{(1)}|^2 \geq \int_{T} |1 + a_\eps +1|^2 + O(\eps^2)\geq C\eps  ,$$
which is clearly impossible. 
Returning to \eqref{eq:first-component-bdry} and \eqref{eq:y-on-rigid-components}, it now follows that $|\tilde{\zeta}^{(i)}_1| = O(\eps)$ and then 
\begin{align*}
   \| \bar{u} - (0, s_\eps) \|_{H^1(\Omega^{(1)}_{\eps})} + \| \bar{u} - (al, t_\eps) \|_{H^1(\Omega^{(2)}_{\eps})} 
  =  O(\sqrt{\eps}),
\end{align*}
where $s_\eps = \frac{1}{\sqrt{\eps}}\zeta^{(1)}_2 $ and $t_\eps = \frac{1}{\sqrt{\eps}}\zeta^{(2)}_2$. \eop 

\begin{lemma}\label{lemma:strong-convergence-special} 
If $a > a_{\rm crit}$ and $\phi = 0$, then there exist sequences $s_{\eps}, t_{\eps} \in \R$ and Lipschitz functions $g_\eps$ as in Lemma \ref{lemma: gamma special case} such that 
$$ \| \bar{u}_{\eps} - (0, s_{\eps}) \|_{H^1(\Omega^{(1)}_{g_\eps})} 
   + \| \bar{u}_{\eps} - (al, t_{\eps}) \|_{H^1(\Omega^{(2)}_{g_\eps})} 
   \to 0. $$    
\end{lemma} 

\Proof Without restriction we only estimate $\bar{y}$ (again dropping subscripts $\eps$) on $\Omega^{(1)}_{g_\eps}$. We claim that we can find a partition $(T_j)_j$, $j=1,\ldots, M_\eps$ of $\Omega^{(1)}_{g_\eps}$ of the form $T_j = \left\{x \in \Omega^{(1)}_{g_\eps}: t_{j-1} \leq x_2 \leq t_{j}\right\}$ for suitable $t_j \in [0,1]$, $j=0,\ldots,M_\eps$ with $t_0=0$ and $t_{M_\eps} =1$ such that the $T_j$ are related through bi-Lipschitzian homeomorphism with uniformly bounded Lipschitz constants to cubes of sidelength $d_j = t_j - t_{j-1} \ge \psi(\eps) \gg \eps$. We will show this at the end of the proof. Recalling that the constant in (\ref{eq:geometric-rigidity}) is invariant under rescaling of the domain and repeating the above arguments in \eqref{eq:y-on-rigid-components} we obtain $R^{(j)} \in SO(2)$ and $\xi^{(j)} \in \R^2$, $j = 1, \ldots, M_{\eps}$, such that
\begin{align*}
  \sum^{M_\eps}_{j=1} \left\|\nabla \bar{y} - R^{(j)}\right\|_{L^2(T_j)}^2 = O(\eps^2)
  \quad\text{and}\quad
  \sum^{M_\eps}_{j=1} d_{j}^{-2}\left\|\bar{y} - R^{(j)} \cdot - \xi^j\right\|_{L^2(T_j)}^2 = O(\eps^2).
\end{align*}
Let $\tilde{T}_j = (t_{j-1},t_j)$ for $j=1,\ldots, M_\eps$ and $T^* = \bigcup^{M_\eps}_{j=1} (t_{j-1} + \frac{d_j}{2}, t_j)$. A standard rescaling argument and the trace theorem yield 
\begin{align}\label{eq: cube traces}
\sum^{M_\eps}_{j=1} d_j^{-1}\left\|\bar{y}(0,\cdot) - R^{(j)} (0,\cdot) - \xi^j\right\|_{L^2(\tilde{T}_j)}^2 = O(\eps^2) .
\end{align}
Similarly as above we calculate the norm in $L^2(T^* \setminus \Gamma^{(1)}_\eps)$ on the trace $\{x_1 = 0\}$ of the piecewise constant function 
\begin{align*}
  \frac{d_j}{2} R^{(j)}_{12}
  = \left( x_1  - R^{(j)}_{11} x_1 - R^{(j)}_{12} \left(\cdot - \frac{d_j}{2}\right) - \xi^{(j)}_1 \right) 
    - \left( x_1  - R^{(j)}_{11} x_1 - R^{(j)}_{12} \cdot - \xi^{(j)}_1 \right) 
\end{align*}
and now find that $\sum^{M_\eps}_{j=1} d_j^2|R^{(i)}_{12}|^2 = O(\eps^2 )$. Consequently, noting that $d_j \ge \psi(\eps) \gg \eps$ for all $j=1,\ldots,M_\eps$ and proceeding as before, we obtain 
\begin{align}\label{eq:NRIdeps} 
  \sum^{M_\eps}_{j=1} d_j^2|R^{(j)} - \Id|^2 = O(\eps^2 ) 
\end{align}
so that 
\begin{align*}
\sum^{M_\eps}_{j=1} \left\|\bar{y} - \id - \xi^j\right\|_{H^1(T_j)}^2 = O(\eps^2). 
\end{align*}
Due to the boundary conditions, \eqref{eq: cube traces} and \eqref{eq:NRIdeps} yield $\sum^{M_\eps}_{j=1} d_j^2 |\xi^j_1|^2 = O(\eps^2)$ and therefore  
\begin{align}\label{eq: xi}
\sum^{M_\eps}_{j=1} \left\|\bar{y} - \id - (0,\xi^j_2)\right\|_{H^1(T_j)}^2 = O(\eps^2). 
\end{align}
We define the stripe $S = (0,\bar{\psi}(\eps)) \times (0,1)$ and note that $S \subset \Omega^{(1)}_{g_\eps}$ by \eqref{eq: components connected}. From Poincar{\'e}'s inequality we obtain a $\zeta \in \R^2$ such that
\begin{align}\label{eq: zeta}
\left\|\bar{y} - \id - \zeta\right\|^2_{H^1(S)} \leq C \left\|\nabla \bar{y} - \Id\right\|^2_{L^2(S)} =  O(\eps^2).
\end{align}
Note that the constant $C$ can be chosen independently of the length of $S$, i.e. independently of $\eps$. Applying \eqref{eq: xi} we may suppose $\zeta=(0,\zeta_2)$. 

Moreover, by \eqref{eq: xi} and \eqref{eq: zeta} there is some $\rho_\eps \in (0, \bar{\psi}(\eps))$ such that the trace on the slice $\Gamma = \left\{\rho_\eps\right\} \times (0,1)$ satisfies
\begin{align*}
\int_\Gamma |\bar{y} - \id - (0,\zeta_2)|^2 = \frac{O(\eps^2)}{\bar{\psi}(\eps)}  
\quad\text{and}\quad
\sum^{M_\eps}_{j=1} \left\|\bar{y} - \id - (0,\xi^j_2)\right\|^2_{L^2(\Gamma \cap \overline{T_j})} = \frac{O(\eps^2)}{\bar{\psi}(\eps)}.  
\end{align*}
We compare the trace on $\Gamma$ and deduce
\begin{align*} 
   \sum^{M_\eps}_{j=1} d_j |\zeta - \xi^j|^2 
   & \leq C\sum^{M_\eps}_{j=1} \Big(\left\|\bar{y} - \id - \xi^j\right\|^2_{L^2(\Gamma \cap \overline{T_j})} + \left\|\bar{y} - \id - \zeta\right\|^2_{L^2(\Gamma \cap \overline{T_j})} \Big)  \\
   &= \frac{O(\eps^2)}{\bar{\psi}(\eps)} = \frac{O(\eps^2)}{\psi(\eps)}.
\end{align*}
Thus, returning to \eqref{eq: xi} we conclude 
\begin{align*}
  \left\|\bar{y} - \id - (0,\zeta_2)\right\|^2_{H^1(\Omega^{(1)}_{g_\eps})} 
   &\leq C\sum^{M_\eps}_{j=1} \left\|\bar{y} - \id - \xi^j\right\|^2_{H^1(T_j)} 
    + C\sum^{M_\eps}_{j=1} d_j^2 |\xi_j - \zeta|^2 \\ 
  & \leq O(\eps^2) + C\sum^{M_\eps}_{j=1} d_j |\xi_j - \zeta|^2 
   = \frac{O(\eps^2)}{\psi(\eps)}
\end{align*}
and finally 
$$\left\|\bar{u} - (0,s_\eps)\right\|^2_{H^1(\Omega^{(1)}_{g_\eps})} = \frac{O(\eps)}{\psi(\eps)} \to 0$$
for $\eps \to 0$, where $s_\eps = \frac{1}{\sqrt{\eps}} \zeta_2$. For $\Omega^{(2)}_{g_\eps}$ we proceed likewise.

To finish the proof it suffices to show the existence of a partition $(T_j)_j$ with the above properties.
Recall that $\Omega^{(1)}_{g_\eps} = \left\{x \in \Omega: 0 < x_1 < g(x_2) -c\eps\right\}$ and $\|g'\|_\infty = \frac{1}{\sqrt{3}}$, $g \ge \psi(\eps)$. Let $r_0 = 0$ and define $r_1, \ldots, r_{M_{\eps}} \in (0, 1)$ inductively by setting $r_{j+1} = r_j + g(r_j)$, so that $r_{M_{\eps}} + g(r_{M_{\eps}}) \ge 1$. Now setting 
$$ T_j = 
   \begin{cases} 
     \{ x \in \Omega^{(1)}_{g_\eps} : r_{j-1} \le x_2 \le r_j \} &\text{for } 1 \le j \le M_{\eps}-1, \\  
     \{ x \in \Omega^{(1)}_{g_\eps} : r_{M_{\eps}-1} \le x_2 \le 1 &\text{for } j = M_{\eps}, 
   \end{cases} $$ 
it is not hard so see that every $T_j$ is related to $\lambda (0,1)^2$ for a suitable $\lambda$ through some bi-Lipschitzian homeomorphism with uniformly bounded Lipschitz constants. By construction, $t_j - t_{j-1} \ge g(t_j) \ge \psi(\eps) \gg \eps$ for $j = 1, \ldots, M_{\eps}$.
\eop

Strong convergence in the subcritical case can be shown along the lines of the proofs of the main linearization results in \cite{Schmidt:08} and \cite{Schmidt:2009}. We include a simplified proof adapted to the present situation here for the sake of completeness. 
\begin{lemma}\label{lemma:strong-convergence-elast} 
If $a < a_{\rm crit}$, then there is a sequence $s_{\eps} \in \R$ such that 
$$ \left\| \bar{u}_{\eps} - (0, s_{\eps}) 
   - F^a \cdot \right\|_{H^1(\Omega_{\rm good})}  
   \to 0. $$ 
where $F^a = \begin{pmatrix} a & 0 \\ 0 & - \frac{a}{3} \end{pmatrix}$. 
\end{lemma} 

\Proof We again drop subscripts $\eps$ if no confusion arises. With the help of the geometric rigidity estimate \eqref{eq:geometric-rigidity} we find by arguing as in the proof of Lemma \ref{lemma:strong-convergence-crack} that 
\begin{align*}
  \| \nabla \bar{y} - R \|_{L^2(\Omega_{\eps})}^2 
  &\le C \int_{\Omega_{\eps} \setminus \bigcup_{\triangle \in \bar{\cal C}_{\eps}} \triangle} 
      W_{\triangle, \chi}(\nabla \tilde{y}) \, dx 
   = O(\eps) 
\end{align*}
for a suitable rotation $R \in SO(2)$ with 
\begin{equation}\label{eq: rot-estimate}
|R \pm \Id| = O(\sqrt{\eps})
\end{equation}
and 
\begin{align*}
  \| \bar{y} \pm \id - \zeta \|_{H^1(\Omega_{\eps})} 
  = O(\sqrt{\eps}) 
\end{align*}
for some $\zeta \in \R^2$ with $\zeta_1 = O(\sqrt{\eps})$ and thus, due to the boundary conditions, 
\begin{align*}
  \| \bar{u} - (0,\zeta_2) \|_{H^1(\Omega_{\eps})} 
  = O(1).  
\end{align*}
In particular, $\bar{u}_{\eps} -(\zeta_{\eps})_{2} \e_2$ converges -- up to passing to a subsequence -- weakly. It now suffices to prove that $\| e(\bar{u}_{\eps}) - F^a\|_{L^2(\Omega_{\eps})} \to 0$, where $e(u) = \frac{(\nabla u)^T + \nabla u}{2}$ denotes the symmetrized gradient, for then the assertion follows from Korn's inequality. 

To this end, we let $V_{\eps}(F) = \frac{1}{\eps} W_{\triangle}(\Id + \sqrt{\eps}F)$ and $V_{\eps,\chi}(F) = V_{\eps}(F) + \frac{1}{\eps}\chi(\Id + \sqrt{\eps}F)$, so that $V_{\eps, \chi}(F) \to \frac{1}{2} D^2 W_{\triangle}(\Id)[F,F] = \frac{1}{2} Q(F)$ uniformly on compact subsets of $\R^{2 \times 2}$. Then by frame indifference (see Lemma \ref{lemma:W-triangle-properties})
\begin{align}\label{eq:W-triangle-V}
\begin{split}
  W_{\triangle,\chi}(\Id + \sqrt{\eps} F) 
  &= W_{\triangle,\chi} \left( \sqrt{(\Id + \sqrt{\eps} F)^T(\Id + \sqrt{\eps} F)} \right) \\ 
  &= \eps V_{\eps,\chi} \left( \frac{F^T + F}{2} + \frac{1}{\sqrt{\eps}} f(\sqrt{\eps} F) \right) 
\end{split}
\end{align}
with $f(F) = \sqrt{(\Id + F)^T (\Id + F)} - \Id - \frac{F^T + F}{2}$, so that $|f(F)| \le C \min\{|F|, |F|^2\}$. Then by Lemma \ref{lemma:quadratic lower bound}(i) and \eqref{eq:W-triangle-V} $V_{\eps,\chi}$ satisfies 
\begin{align}\label{eq: dist-V_eps}
\begin{split}
  V_{\eps,\chi}\left( \frac{F^T + F}{2} + \frac{1}{\sqrt{\eps}} f(\sqrt{\eps} F) \right) 
  &\ge \frac{c}{\eps} \dist^2(\Id + \sqrt{\eps} F, O(2)) + \frac{1}{\eps}\chi(\Id + \sqrt{\eps}F)\\ 
  &\ge \frac{c}{\eps} \dist^2(\Id + \sqrt{\eps} F, SO(2))\\
  &\ge \frac{c}{\eps} \left|\sqrt{(\Id + \sqrt{\eps} F)^T(\Id + \sqrt{\eps} F)} - \Id \right|^2 \\ 
  &= c \left| \frac{F^T + F}{2} + \frac{1}{\sqrt{\eps}} f(\sqrt{\eps} F) \right|^2.  
\end{split} 
\end{align}
In the sequel we set $A_{\eps}(F) = \frac{F^T + F}{2} + \frac{1}{\sqrt{\eps}} f(\sqrt{\eps} F)$. Choose convex functions $\psi_k : \R^{2 \times 2} \to \R$ with linear growth at infinity such that $\psi_1 \le \psi_2 \le \ldots$ and $\psi_k(F) \to \frac{1}{2} Q(F)$ uniformly on compact subsets of $\R^{2 \times 2}$. The previous quadratic estimate on $V_{\eps,\chi}(A_{\eps}(F))$ from below and the fact that $V_{\eps, \chi} \to \frac{1}{2} Q$ uniformly on compacts then shows that we can also choose $\delta > 0$ and a sequence $r_k \to \infty$ such that 
\begin{align*}
  V_{\eps,\chi} \left( A_{\eps}(F) \right) 
  - \delta \chi_{\{|A_{\eps}(F)| \ge r_k\}} |A_{\eps}(F)|^2 
  \ge \psi_k\left( A_{\eps}(F) \right) - \frac{1}{k},
\end{align*}
whenever $\eps$ (depending on $k$) is sufficiently small.  

With \eqref{eq:W-triangle-V} we now obtain that
\begin{align*}
  \frac{1}{\eps} \int_{\Omega_{\eps}} W_{\triangle,\chi}( \bar{y} ) \, dx 
  &= \int_{\Omega_{\eps}} 
     V_{\eps,\chi}\left( A_{\eps}(\nabla \bar{u}) \right) \, dx \\ 
  &\ge \int_{\Omega_{\eps}} 
     \psi_k \left( A_{\eps}(\nabla \bar{u}) \right) \, dx 
     + \delta \int_{\Omega_{\eps}} \chi_{\{ |A_{\eps}(\nabla \bar{u})| \ge r_k\}} 
     |A_{\eps}(\nabla \bar{u})|^2 \, dx 
     - \frac{1}{k}. 
\end{align*}
As $\psi_k$ has linear growth at infinity and $\frac{1}{\sqrt{\eps}} f(\sqrt{\eps} \nabla \bar{u}_{\eps}) \le C \min\{|\nabla \bar{u}_{\eps}|, \sqrt{\eps}|\nabla \bar{u}_{\eps}|^2\}$, $\nabla \bar{u}_{\eps}$ bounded in $L^2$, by splitting the integration into two parts according to $|\nabla \bar{u}_{\eps}| \leq M$ or $|\nabla \bar{u}_{\eps}| > M$ and eventually sending $M$ to infinity, we find 
\begin{align*}
 \liminf_{\eps \to 0} \int_{\Omega_{\eps}} 
     \psi_k \left( A_{\eps}(\nabla \bar{u}_{\eps}) \right) \, dx = \liminf_{\eps \to 0}  \int_{\Omega_{\eps}} \psi_k \left( e(\bar{u}_{\eps}) \right) \, dx.
\end{align*}
When $\bar{u}_{\eps} - (\zeta_{\eps})_2\e_2\weakly u$ in $H^1$, by Theorem \ref{theo:limiting-energy} it then follows that 
\begin{align*}
  \frac{\alpha l a^2}{\sqrt{3}} 
  &= \lim_{\eps \to 0} \frac{4}{\sqrt{3}} \int_{\Omega_{\eps}} V_{\eps,\chi}\left( A_{\eps}(\nabla \bar{u}_{\eps}) \right) \, dx \\ 
  &\ge \liminf_{\eps \to 0} \frac{4}{\sqrt{3}} \int_{\Omega} \chi_{\{ \dist(x, \partial \Omega) \ge k^{-1}\}} \psi_k \left( e(\bar{u}_{\eps}) \right) \, dx \\ 
  &\qquad  + \limsup_{\eps \to 0} \frac{4 \delta}{\sqrt{3}} \int_{\Omega_{\eps}} \chi_{\{|A_{\eps}(\nabla \bar{u}_{\eps})| \ge r_k\}} 
     |A_{\eps}(\nabla \bar{u}_{\eps})|^2 \, dx 
     - \frac{4}{\sqrt{3}k}. 
\end{align*}
Using that by convexity of $\psi_k$ the first term on the right hand side is lower semicontinuous in $\nabla \bar{u}_{\eps}$ and that $\chi_{\{ \dist(\cdot, \partial \Omega) \ge k^{-1}\}} \psi_k \to \frac{1}{2} Q$ monotonically, we finally find by letting $k \to \infty$ 
\begin{align}\label{eq: energy of u}
\begin{split}
  \frac{\alpha l a^2}{\sqrt{3}} 
  &\ge \frac{2}{\sqrt{3}} \int_{\Omega} Q\left( e(u) \right) \\
      &\qquad  + \lim_{k \to \infty} \limsup_{\eps \to 0} \frac{4 \delta}{\sqrt{3}} \int_{\Omega_{\eps}} \chi_{\{|A_{\eps}(\nabla \bar{u}_{\eps})| \ge r_k\}} 
    |A_{\eps}(\nabla \bar{u}_{\eps})|^2 \, dx. 
  \end{split} 
\end{align} 
A slicing and convexity argument similar to (\ref{eq: elastic energy estimate}) now shows that $\frac{2}{\sqrt{3}} \int_{\Omega} Q(e(w)) \ge \frac{\alpha l a^2}{\sqrt{3}}$ for all $w \in H^1$ subject to $w_1(0, x_2) = 0$ and $w_1(l, x_2) = al$ and thus 
\begin{align*}
  \lim_{k \to \infty} \limsup_{\eps \to 0} \frac{4 \delta}{\sqrt{3}} \int_{\Omega_{\eps}} \chi_{\{|A_{\eps}(\nabla \bar{u}_{\eps})| \ge r_k\}} 
    |A_{\eps}(\nabla \bar{u}_{\eps})|^2 \, dx = 0, 
\end{align*} 
or, in other words, $|A_{\eps}(\nabla \bar{u}_{\eps})|^2$ is equiintegrable. By the estimate $|V_{\eps, \chi}(F)| = |\frac{1}{\eps}W_{\triangle, \chi}(\Id + \sqrt{\eps} F)| \leq C(1 + |F|^2)$, (\ref{eq: dist-V_eps}) shows that also 
$$ \frac{c}{\eps} \dist^2(\nabla \bar{y}_{\eps}, SO(2)) \leq V_{\eps, \chi}(A_{\eps}(\nabla \bar{u}_{\eps})) $$
is equiintegrable, so that by the discussion following Equation \eqref{eq:geometric-rigidity} in fact we may assume that $\frac{1}{\eps} \| \nabla \bar{y}_{\eps} - R \|_{L^2(\Omega_{\eps})}^2$ is equiintegrable, too, and $|R - \Id| = O(\sqrt{\eps})$  by (\ref{eq: rot-estimate}). But then also $|\nabla \bar{u}_{\eps}|^2$ is equiintegrable and this together with (\ref{eq: energy of u}) yields 
\begin{align*}
\lim_{\eps \to 0}\frac{2}{\sqrt{3}} \int_{{\Omega_{\eps}}} Q(e(\bar{u}_{\eps})) = \frac{2}{\sqrt{3}} \int_{{\Omega}} Q(e(u)) = \frac{\alpha l a^2}{\sqrt{3}}.
\end{align*}
For some $\delta > 0$ small enough we finally obtain that 
\begin{align*}
  \frac{\alpha l a^2}{\sqrt{3}} 
  &= \frac{2}{\sqrt{3}} \int_{\Omega} Q(F^a) \, dx \\ 
  &= \inf \bigg\{ \frac{2}{\sqrt{3}} \int_{\Omega} Q(e(w)) - \delta |e(w) - F^a|^2 \, dx : \\ 
  &\qquad \qquad \qquad w \in H^1(\Omega), w(0,x_2) = 0, w(l, x_2) = al \bigg\} \\ 
  &\le \liminf_{\eps \to 0} \frac{2}{\sqrt{3}} \int_{{\Omega_{\eps}}} Q(e(\bar{u}_{\eps})) - \delta | e(\bar{u}_{\eps}) - F^a |^2 \, dx \\ 
  &= \frac{\alpha l a^2}{\sqrt{3}} - \delta \limsup_{\eps \to 0} \|e(\bar{u}_{\eps}) - F^a\|^{2}_{L^2{ (\Omega_{\eps})}}
  \end{align*} 
and therefore $\lim_{\eps \to 0} \|e(\bar{u}_{\eps}) - F^a\|^{2}_{L^2{ (\Omega_{\eps})}} = 0$ indeed. \eop 

After all these preparatory lemmas, the proof of our main limiting result Theorem \ref{theo: conv-of-minimizers} is now straightforward. \smallskip

\noindent {\em Proof of Theorem \ref{theo: conv-of-minimizers}.} 
Choose $s_{\eps}$ as in Lemmas \ref{lemma:strong-convergence-elast} if $a < a_{\rm crit}$, $p_{\eps}$, $s_{\eps}$ and $t_{\eps}$ as in \eqref{eq:p-eps} and Lemma \ref{lemma:strong-convergence-crack} if $a > a_{\rm crit}$ and $\phi \ne 0$ and finally $g_\eps$ and $s_\eps$ and $t_\eps$ as in Lemma \ref{lemma:strong-convergence-special} if $a > a_{\rm crit}$ and $\phi = 0$. By Lemmas \ref{lemma:strong-convergence-elast}, \ref{lemma:strong-convergence-crack} and \ref{lemma:strong-convergence-special}, $\bar{u}_{\eps}$ can be extended as an $H^1$-function from $\Omega_{\eps}$ to $\Omega$, $\Omega^{(i)}_{\eps}$ to $\Omega^{(i)}$, $i = 1,2$, or $\Omega^{(i)}_{g_\eps}$ to $\Omega^{(i)}[g_\eps]$, $i=1,2$, respectively, such that still, respectively, 
\begin{align}
  \left\| \bar{u}_{\eps} - (0, s_{\eps}) - F^a \cdot \right\|_{H^1(\Omega)}  
  &\to 0, \label{eq: convergence1} \\ 
  \| \bar{u}_{\eps} - (0, s_{\eps}) \|_{H^1(\Omega^{(1)})} + \| \bar{u}_{\eps} - (al, t_{\eps}) \|_{H^1(\Omega^{(2)})} 
  &\to 0, \label{eq: convergence2} \\ 
  \| \bar{u}_{\eps} - (0, s_{\eps}) \|_{H^1(\Omega^{(1)}[g_\eps])} + \| \bar{u}_{\eps} - (al, t_{\eps}) \|_{H^1(\Omega^{(2)}[g_\eps])} 
  &\to 0. \label{eq: convergence3}
\end{align}
This completes the proof as by Lemma \ref{lemma:triangle-healing} we also still have $|\{x \in \Omega_{\eps}: \bar{u}_{\eps}(x) \ne \tilde{u}_{\eps}(x)\}| = O(\eps)$. \eop

Finally, we give the proof of Corollary \ref{cor: limiting-deformations}. 

\noindent {\em Proof of Corollary \ref{cor: limiting-deformations} .} 
First, let $(y_\eps)$ be a minimizing sequence satisfying \eqref{eq:almost-min}. Then by Theorem \ref{theo: conv-of-minimizers} we obtain \eqref{eq: convergence1}, \eqref{eq: convergence2} or \eqref{eq: convergence3}, respectively. Taking the condition $\sup_\eps \left\|y_\eps\right\|_\infty < \infty$ into account, in the cases (i) and (ii) we get $\sup_\eps |s_\eps| < \infty$ and $\sup_\eps |s_\eps|, \sup_\eps |t_\eps| < \infty$ such that, passing to subsequences, we obtain $s_\eps \to s$ and $s_\eps \to s, t_\eps \to t, p_\eps \to p$, respectively, for suitable constants $s,t \in \R$, $p \in (0,l)$. In (iii) we first note that up to subsequences $g_\eps$ converges uniformly to some Lipschitz function $g:(0,1) \to [0,l]$ satisfying $|g'|\leq \frac{1}{\sqrt{3}}$ a.e. Then using again the uniform bound $\sup_\eps \left\|y_\eps\right\|_\infty < \infty$ we get constants $s,t$ such that ${s}_\eps \to s$ and $t_\eps \to t$ up to subsequences. It follows that $\tilde{u}_{\eps} \to u$ as given in (i), (ii) and (iii), respectively. 

Conversely, we assume that $u$ is given as in Corollary \ref{cor: limiting-deformations} and show that there is a minimizing sequence $(y_\eps)$ satisfying \eqref{eq:almost-min} with $\tilde{u}_\eps \to u$ in measure. For (i) and (ii) this is obvious by the proof of Theorem \ref{theo:limiting-energy} taking the configurations in \eqref{eq: elastic minimizer} and \eqref{eq: crack minimizer} up to suitable translations. For given $u$ in (iii) with corresponding function $g$ and constants $s,t$ we approximate $g:(0,1) \to [0,l]$ uniformly by Lipschitz functions $g_\eps:(0,1)\to(0,l)$ being affine on intervals of length $\frac{\sqrt{3}\eps}{2}$ with $g'_\eps = \pm \frac{1}{\sqrt{3}}$ a.e. 
We set 
$$y_\eps(x) = \begin{cases} x + (0, \sqrt{\eps} s), & \mbox{if } 0 < x_1 < g_\eps(x_2), \\ 
                x + (a_{\eps} l, \sqrt{\eps} t), & \mbox{if } g_\eps(x_2) < x_1 < l, 
\end{cases}
$$
so that $\tilde{u}_\eps = \frac{y_\eps - \id}{\sqrt{\eps}} \to u$ in measure. As in the proof of Theorem \ref{theo:limiting-energy}, except for negligible contributions of the boundary layers, ${\cal E}^\chi_\eps (y_\eps)$ is given by the energy of the springs intersected transversally by $\text{graph}(g_\eps)$. These springs are elongated by a factor scaling with $\frac{1}{\sqrt{\eps}}$ yielding a contribution $\eps \beta$ in the limit. It is elementary to see that on every stripe in $\e_1$ direction of length $\frac{\sqrt{3}\eps}{2}$ the graph intersects two springs, and consequently ${\cal E}^\chi_\eps (y_\eps) \to \frac{4\beta}{\sqrt{3}}$. \eop


 \typeout{References}

\end{document}